\documentclass[a4paper,11pt,oneside]{article}
\usepackage{amsmath,amssymb,amsthm}
\usepackage[frame,cmtip,curve,arrow,matrix,line,graph]{xy}
\usepackage[dvips]{graphicx}

\newtheorem{thm}{Theorem}[section]
\newtheorem{mainthm}[thm]{Main Theorem}
\newtheorem{cor}[thm]{Corollary}
\newtheorem{prop}[thm]{Proposition}
\newtheorem{defn}[thm]{Definition}
\newtheorem{lem}[thm]{Lemma}

\def\Stab{\mathop{\mathrm{Stab}}\nolimits}

\def\NE{\mathop{\mathrm{NE}}\nolimits}

\def\Coh{\mathop{\mathrm{Coh}}\nolimits}
\def\Spec{\mathop{\mathrm{Spec}}\nolimits}
\def\End{\mathop{\mathrm{End}}\nolimits}

\def\D{D^b}
\def\H{\mathcal{H}}

\def\T{\mathcal{T}}
\def\C{\mathbb{C}}
\def\R{\mathbb{R}}

\def\e{\varepsilon}
\newcommand{\N}{\mathcal{N}}
\newcommand{\M}{\mathcal{M}}
\def\min{\mathop{\mathrm{min}}\nolimits}

\def\fr{\mathop{\mathrm{fr}}\nolimits}
\def\tor{\mathop{\mathrm{tor}}\nolimits}

\def\Hom{\mathop{\mathrm{Hom}}\nolimits}
\newcommand{\K}{{{K}}}
\newcommand{\p}{\mathcal{P}}
\newcommand{\mo}{\mathcal{O}}
\newcommand{\E}{\mathcal{E}}

\newcommand{\PP}{\mathbb{P}}
\newcommand{\Q}{\mathbb{Q}}

\newcommand{\A}{\mathcal{A}}
\newcommand{\Z}{\mathbb{Z}}
\newcommand{\Pic}{\operatorname{Pic}}
\newcommand{\GL}{\operatorname{GL}}
\newcommand{\grp}{{\widetilde{\GL}^+}(2,\R)}
\def\ns{\mathop{\mathrm{NS}}\nolimits}
\def\Im{\mathop{\mathrm{Im}}\nolimits}
\def\Re{\mathop{\mathrm{Re}}\nolimits}
\def\im{\mathop{\mathrm{im}}\nolimits}

\def\ker{\mathop{\mathrm{ker}}\nolimits}
\def\coker{\mathop{\mathrm{coker}}\nolimits}
\def\dim{\mathop{\mathrm{dim}}\nolimits}
\def\dimv{\mathop{\mathrm{\underline{dim}}}\nolimits}

\def\rk{\mathop{\mathrm{rk}}\nolimits}

\def\ch{\mathop{\mathrm{ch}}\nolimits}
\def\Mod{\mathop{\mathrm{mod}}\nolimits}
\def\Amp{\mathop{\mathrm{Amp}}\nolimits}
\def\Vec{\mathop{\mathrm{Vec}}\nolimits}
\def\id{\mathop{\mathrm{id}}\nolimits}

\def\exp{\mathop{\mathrm{exp}}\nolimits}
\begin{document}
\title{Moduli of Bridgeland semistable objects on $\PP^2$}
\author{Ryo Ohkawa}
\date{}

\maketitle

\section{Introduction}
Let $X$ be a smooth projective 
surface and $\D(X)$ the bounded 
derived category of coherent 
sheaves on $X$. We study 
Bridgeland stability conditions 
$\sigma$ on $\D(X)$. We show that if
a stability condition $\sigma$
has a certain property, 
the moduli space of 
$\sigma$-(semi)stable 
objects in $\D(X)$ coincides with 
a certain moduli space of 
Gieseker-(semi)stable 
coherent sheaves on $X$. On the other 
hand, when $X$ has a full strong
exceptional collection, we define
the notion of $\sigma$ being 
''algebraic'', and we show that
for any algebraic stability 
condition $\sigma_{\text{alg}}$, 
the moduli space of $\sigma_{
\text{alg}}$-(semi)stable objects in 
$\D(X)$ coincides with a certain moduli 
space of modules over a finite dimensional 
$\C$-algebra. Using these 
observations, we 
construct moduli spaces of 
Gieseker-(semi)stable 
coherent sheaves on $\PP^2$ as moduli 
spaces of certain
modules (Theorem~\ref{our}).
This gives a new proof (\S~\ref{comp})
of Le Potier's result \cite{P} and 
establishes some related results 
(\S~\ref{cwc}).

\subsection{Bridgeland stability 
conditions}
The notion of stability conditions 
on a triangulated category $\T$ 
was introduced in \cite{B2} to give 
the mathematical framework for the 
Douglas's work on $\Pi$-stability.
Roughly speaking, 
it consists of data $\sigma=
(Z,\A)$, where $Z$ is a group 
homomorphism from the Grothendieck 
group $K(\T)$ to the complex number 
field $\C$, $\A$ is a full 
abelian subcategory of $\T$ and 
these data should have some 
properties (see Definition~\ref{stc}). 
Then Bridgeland \cite{B2} showed 
that the set of some good stability 
conditions has a structure of a 
complex manifold.  
This set is denoted by 
$\Stab(X)$ when $\T=\D(X)$.
An element $\sigma$ of $\Stab(X)$ 
is called a Bridgeland stability 
condition on $X$. 
For a full abelian 
subcategory $\A\subset\T$, $\Stab(\A)$ 
denotes the subset of $\Stab(X)$ 
consisting of all stability conditions 
of the form $\sigma=(Z,\A)$.

Let $K(X)$ be the Grothendieck group of 
$X$. For $\alpha\in K(X)$, the Chern character 
of $\alpha$ is the element 
$\ch(\alpha):=(\rk(\alpha),c_1(\alpha),\ch_2(\alpha))$
of the lattice
$\N(X):=\Z\oplus\ns(X)\oplus\frac{1}{2}\Z$.
For $\sigma=
(Z,\A)\in\Stab(X)$, we consider the moduli 
functor $\M_{\D(X)}(\ch(\alpha),\sigma)$ of 
$\sigma$-(semi)stable objects $E$ in $\A$ 
with $\ch(E)=\ch(\alpha)$.

\subsection{Geometric Bridgeland 
stability conditions}
For $\beta, \omega\in\ns(X)
\otimes\R$ such that $\omega$ is 
in the ample cone $\Amp(X)$, 
we consider a pair
$\sigma_{(\beta,\omega)}=
(Z_{(\beta,\omega)},\A_{(
\beta,\omega)})$ as in 
\cite{ABL}, where 
$Z_{(\beta,\omega)}\colon
K(X)\to\C$ is a group 
homomorphism and $\A_{(\beta,
\omega)}$ is a full abelian
subcategory of $\D(X)$
defined from $\beta$ and 
$\omega$
(see Definition~\ref{def of 
geometric} for details). 
It is shown in \cite{ABL}
that $\sigma_{(\beta,\omega)}$
is a Bridgeland stability 
condition if $\beta,\omega
\in\ns(X)\otimes\Q$. 
For general $\beta,\omega\in
\ns(X)\otimes\R$, we do not
know whether $\sigma_{(\beta,\omega)}$
belongs to $\Stab(X)$ or not
(cf. \S~\ref{ABL}).

Let $\grp$ be the universal cover
of the group $\GL^+(2,\R):=
\{T\in\GL(2,\R)\mid\det T>0\}$. 
The group $\grp$ acts on $\Stab(X)$ 
in a natural way (cf. \S~\ref{grp}).
Two stability conditions 
$\sigma$ and $\sigma'$ 
are said to be $\grp$-equivalent
if $\sigma$ and $\sigma'$ are
in a single orbit of this
action. 
In such cases $\sigma$ and
$\sigma'$ correspond to 
isomorphic moduli functors of
semistable objects.
$\sigma\in\Stab(X)$ is 
said to be geometric if 
$\sigma$ is 
$\grp$-equivalent to 
$\sigma_{(\beta,\omega)}$
for some $\beta,\omega
\in\ns(X)\otimes\R$ with 
$\omega\in\Amp(X)$.
We have a criterion due to
\cite{B3} for 
$\sigma\in\Stab(X)$ to be 
geometric (Proposition~\ref
{geometric1}).

On the other hand, 
for an integral ample divisor 
$\omega$ and $\beta\in\ns(X)\otimes
\Q$, we consider
$(\beta,\omega)$-twisted
Gieseker-stability 
of torsion free sheaves on $X$,
which was introduced in 
\cite{MW} generalizing 
the Gieseker-stability. 
For $\alpha\in K(X)$, 
we assume $\rk(\alpha)>0$
and consider the moduli 
functor $\M_X(\ch(\alpha),\beta,
\omega)$ of $(\beta,
\omega)$-semistable 
sheaves $E$ with 
$\ch(E)=\ch(\alpha)$. 
There is a scheme $M_X(\ch(\alpha),\beta,
\omega)$ which 
corepresents $\M_X(\ch(\alpha),
\beta,\omega)$ \cite{MW}, 
and is called the moduli
space 
(cf. Definition~\ref{corep}). 

One of our main results is the 
following.
\begin{thm}
\label{geomintro}
Let $\omega$ be an
integral ample divisor,
$\beta\in\ns(X)
\otimes\Q$ and $\alpha\in 
K(X)$ with $\rk(\alpha)>0$. 
Take a real number $t$ with
$0<t\le 1$ and assume that
$\sigma_{(\beta,t\omega)}\in
\Stab(X)$. If  
$0<c_1(\alpha)\cdot
\omega-\rk(\alpha)\beta\cdot
\omega\le\min\{t, \frac{1}{\rk
(\alpha)}\}$
then the moduli space 
$M_X(\ch(\alpha),\beta-
\frac{1}{2}K_X,\omega)$ 
corepresents the moduli functor 
$\M_{\D(X)}(\ch(\alpha),
\sigma_{(\beta,t\omega)})$.
\end{thm}
A proof of Theorem~\ref{geomintro}
will be given in \S~\ref{moco}.
Similar results are obtained 
by \cite{B3} and \cite{To} 
when $X$ is a K3 surface,
but our choices of $\omega$
and $\beta$ are different 
from theirs.

\subsection{Algebraic Bridgeland 
stability conditions}
\label{inalgbri}
For a finite dimensional $\C$-algebra 
$B$, mod-$B$ denotes the abelian 
category of finitely generated right 
$B$-modules and $K(B)$ denotes the 
Grothendieck group. For any $B$-module
$N$, we denote by $[N]$
the image of $N$ by the map mod-$B\to
K(B)$. 
King \cite{Ki1} introduced the 
notion of $\theta_B$-stability 
of $B$-modules, where $\theta_B$ is 
a group homomorphism 
$\theta_B\colon K(B)\to\R$.
It is shown in \cite{Ki1} that 
the moduli space 
$M_{B}(\alpha_B,\theta_B)$ of 
$\theta_B$-semistable $B$-modules $N$ 
with $[N]=\alpha_B$ exists,
for any $\alpha_B\in K(B)$ and 
$\theta_B\in\alpha_B^\perp:=\{\theta_B\in
\Hom_{\Z}(K(B),\R)\mid \theta_B
(\alpha_B)=0\}$. 

When $X$ has a full strong 
exceptional collection 
$\mathfrak{E}=\left(E_0,\ldots,E_n
\right)$ in
$\D(X)$ (cf. \S~\ref{algbri}), 
we put $\E=\oplus_iE_i$ and consider the 
finite dimensional $\C$-algebra 
$B_{\E}=\End_X(\E)$.
Then by Bondal's Theorem \cite{Bo}, 
the functor $\mathbf R\Hom_{X}(\E,\ \cdot\ )$
gives an equivalence of triangulated 
categories $\Phi_{\E}\colon\D(X)\cong
\D(B_{\E})$, where $\D(B_{\E})$ is the 
bounded derived category of mod-$B_{\E}$. 
$\Phi_{\E}$ induces an isomorphism of 
the Grothendieck groups 
$\varphi_{\E}\colon K(X)\cong K(B_{\E})$. 
Let $\A_{\E}$ be the full abelian
subcategory of $\D(X)$ corresponding to
mod-$B_{\E}\subset\D(B_\E)$ by 
$\Phi_{\E}$. $\sigma\in\Stab(X)$ 
is called an algebraic Bridgeland 
stability condition associated 
to $\mathfrak{E}=\left(E_0,\ldots,E_n
\right)$
if $\sigma$ is $\grp$-equivalent to 
$(Z,\A_\E)$ for some $Z\colon
K(X)\to \C$.

For any $\sigma=(Z,\A_{\E})\in\Stab
(\A_{\E})$ and $\alpha\in K(X)$, 
we associate the group 
homomorphism 
$\theta^{\alpha}_Z\colon K(B_\E)
\to\R$ defined by 
$$\theta^{\alpha}_Z(\beta)=
\begin{vmatrix}
\Re Z(\varphi_\E^{-1}(\beta)) & \Re Z(\alpha)\\
\Im Z(\varphi_\E^{-1}(\beta)) & \Im Z(\alpha)
\end{vmatrix}$$
for $\beta\in K(B_\E)$. 
Clearly $\theta^\alpha_Z
\in\varphi_\E(\alpha)^\perp$,
so we have the moduli space
$M_{B_\E}(\varphi_\E(\alpha),
\theta^\alpha_Z)$.
\begin{prop}\label{algintro}
The moduli space 
$M_{B_{\E}}(\varphi_\E(\alpha),
\theta^\alpha_Z)$ of $B_{\E}$-modules
corepresents the 
moduli functor $\M_{\D(X)}(\ch(\alpha),
\sigma)$ for any $\alpha\in K(X)$ and 
$\sigma=(Z,\A_{\E})\in\Stab
(\A_{\E})$.
\end{prop}
A proof of Proposition~\ref{algintro}
will be given in \S~\ref{algbri}.

\subsection{Application in the 
case $X=\PP^2$} 
We prove that there exist 
Bridgeland stability 
conditions on $\PP^2$ which 
are both geometric and algebraic 
by using the criterion 
Proposition~\ref{geometric1}.

The Neron-Severi group $\ns(\PP^2)$ 
of $\PP^2$ is generated by the 
hyperplane class $H$. Hence 
when $X=\PP^2$ the
twisted Gieseker-stability 
coincides with the 
classical one defined by $H$. 
We sometimes identify $\ns(\PP^2)$ with $\Z$
by the map $\beta\mapsto \beta\cdot H$.
For $\alpha\in K(\PP^2)$ with
$\rk(\alpha)>0$, we consider
the moduli space $M_{\PP^2}(\ch(\alpha),H)$
and $\sigma_{(bH,tH)}$ for 
$b,t>0$.

On the other hand, 
for each $k\in\Z$ there exist 
full strong exceptional 
collections on $\PP^2$
\begin{gather*}
\mathfrak{E}_k:=\left(\mo_{\PP^2}(k+1), 
\Omega_{\PP^2}^1(k+3), 
\mo_{\PP^2}(k+2)\right)\ \ 
\text{ and } \ \ 
\mathfrak{E}'_k:=\left(\mo_{\PP^2}(k),
\mo_{\PP^2}(k+1),\mo_{\PP^2}(k+2)\right).
\end{gather*}
We put 
$\E_k:=\mo_{\PP^2}(k+1)\oplus
\Omega_{\PP^2}^1(k+3)\oplus
\mo_{\PP^2}(k+2)$ and 
$\mathcal{E}'_k
:=\mo_{\PP^2}(k)\oplus
\mo_{\PP^2}(k+1)\oplus\mo_{\PP^2}(k+2)$.
Up to natural isomorphism, 
$\End_{\PP^2}(\E_k)$ and $\End_{\PP^2}(\E'_k)$
do not depend on $k$, hence we identify and denote
them by $B$ and $B'$ respectively.
Using the notation in \S~\ref{inalgbri}, we put 
\begin{equation*}
\Phi_k:=\Phi_{\E_k}
\colon\D(\PP^2)
\cong\D(B),\hspace{5mm}
\Phi'_k:=\Phi_{\mathcal{E}'_k}\colon\D(\PP^2)
\cong\D(B'),
\end{equation*}
induced isomorphisms
$\varphi_k:=
\varphi_{\E_k}\colon 
K(\PP^2)\cong K(B)$,
${\varphi'_k}:=\varphi_{\mathcal{E}'_k}
\colon K(\PP^2)\cong K(B')$ and
hearts of induced bounded t-structures
$\A_k:=\A_{\E_k}\subset
\D(\PP^2)$,  
$\A'_k:=\A_{\mathcal{E}'_k}\subset
\D(\PP^2)$.  

For $\alpha\in K(\PP^2)$ and 
$\theta\in\alpha^{\perp}:=
\big\{\theta\in\Hom_{\Z}(K(\PP^2),\R)
\mid\theta(\alpha)=0\big\}$,
we put $$ 
\theta_k:=\theta\circ
\varphi^{-1}_k\colon K(B)
\to\R, \ \  
\theta'_k:=\theta
\circ{\varphi'_k}^{-1}\colon K(B')
\to\R.$$ There exists  
$\theta\in\alpha^{\perp}$ 
such that 
$\Phi'_1\circ\Phi^{-1}_0$ and $\Phi_1
\circ{\Phi'_1}^{-1}$ induce the 
following isomorphisms
(Proposition~\ref{tilt})
\begin{equation}\label{nim}
M_{B}(-\varphi_0(\alpha),\theta_0)
\cong M_{B'}(-\varphi'_1(\alpha),
\theta'_1)\cong M_{B}(-\varphi_1(\alpha),
\theta_1).
\end{equation}

We find algebraic Bridgeland 
stability conditions 
$\sigma^b=(Z^b,\A_1)\in
\Stab(\A_1)$ parametrized 
by real numbers $b$ with 
$0<b<1$ such that for each 
$b$ there exist an element 
$g\in\grp$ and $t>0$ satisfying
\begin{equation}\label{ga}
\sigma^bg=\sigma_{(bH,tH)},
\end{equation}
where $g$ and $t>0$ may depend
on $b$.
Then $M_{B}(-\varphi_1(\alpha),
\theta^{\alpha}_{Z^b})$
corepresents the moduli functors 
$\M_{\D(\PP^2)}(-\ch(\alpha),\sigma^b)$
by Proposition~\ref{algintro}.
Furthermore by (\ref{ga}) and
Theorem~\ref{geomintro},
$M_{\PP^2}(\ch(\alpha),H)$ also
corepresents the same moduli
functor for suitable choice of
$b$.
From these facts and  
isomorphisms (\ref{nim}), we
have our main results
(see \S~\ref{main} for the
choice of $\theta\in\alpha^{\perp}$).
We denote by $\ \cdot\ [1]$
the shift functor
$\D(\PP^2)\to\D(\PP^2)
\colon E\mapsto E[1]$.
\begin{mainthm}
\label{inmain} 
For $\alpha\in K(\PP^2)$ 
with $c_1(\alpha)=sH$, 
assume $0<s\le \rk(\alpha)$ and 
$\ch_2(\alpha)<
\frac{1}{2}$. Then there exists 
$\theta\in\alpha^\perp$
such that $\Phi_1(\ \cdot\ [1])$, 
${\Phi'_1}(\ \cdot\ [1])$ 
and $\Phi_0(\ \cdot\ [1])$ induce the 
following isomorphisms.\\ 
\emph{(i)} $M_{\PP^2}(\ch(\alpha),H)
\cong M_{B}(-\varphi_1(\alpha),
\theta_1)\colon E
\mapsto \Phi_1(E[1])$\\
\emph{(ii)} $M_{\PP^2}(\ch(\alpha),H)
\cong M_{B'}(-\varphi'_1(\alpha),
\theta'_1)\colon E
\mapsto \Phi'_1(E[1])$\\
\emph{(iii)} $M_{\PP^2}(\ch(\alpha),H)
\cong M_{B}(-\varphi_0(\alpha),
\theta_0)\colon E
\mapsto \Phi_0(E[1])$.\\
These isomorphisms keep 
open subsets 
consisting of stable objects. 
\end{mainthm}
We remark that if we assume
$0<s\le\rk(\alpha)$ and  
$M_{\PP^2}(\ch(\alpha),H)\neq\emptyset$
in Main Theorem~\ref{inmain}, 
then we have $$\dim M_{\PP^2}(\ch(\alpha),H)=s^2
-\rk(\alpha)^2+1-2\rk(\alpha)\ch_2(\alpha)\ge 0.$$
Hence we have $\ch_2(\alpha)\le\frac{1}{2}$, 
and $\ch_2(\alpha)=\frac{1}{2}$ if and 
only if $M_{\PP^2}(\ch(\alpha),H)=\{\mo_
{\PP^2}(1)\}$. In this case, similar isomorphisms
hold via $\Phi_1(\ \cdot\ [1])$ in (i),
$\Phi'_1$ in (ii) and $\Phi_0$ in (iii)
respectively. A proof of Main Theorem~\ref{inmain}
will be given in \S~\ref{quiver}.

(ii) is obtained by Le Potier 
\cite{P} (cf. \cite[\S~4]{KW} 
and \cite[Theorem~14.7.1]{P2})
by a different method. 



\subsection{Wall-crossing phenomena}

In \S~\ref{cwc} we consider
the case $\rk(\alpha)
=1$, $c_1(\alpha)=H$ and 
$\ch_2(\alpha)=\frac{1}{2}-n$
with $n\ge 1$.
By Main Theorem~\ref{inmain} we have
$$M_{\PP^2}(\ch(\alpha),H)
\cong M_{B}(-\varphi_0(\alpha),
\theta_0)\cong M_{B}(
-\varphi_1(\alpha),\theta_1)$$
for some $\theta\in\alpha^\perp$. 
We study how $M_{B}(-\varphi_k(\alpha),
\theta_k^{\dagger})$ changes
when $\theta^{\dagger}_k\in\varphi_k(\alpha)^{\perp}$
varies for $k=0,1$,
where $\varphi_k(\alpha)^{\perp}:=\{
\theta_k\in\Hom_\Z(K(B),\R)
\mid \theta_k(\varphi_k(\alpha))
=0\}$.  
We define a wall-and-chamber
structure on $\varphi_k(\alpha)^{\perp}$
as follows (cf. \S~\ref{main}). 
Within $\varphi_k(\alpha)^{\perp}$, 
there are finitely many
rays corresponding to certain
$B$-modules. In our case, a
ray may be called a wall, since
$\varphi_k(\alpha)^{\perp}\cong\R^2$.
Let $W_k$ be the union of such 
rays. A connected component
of the complement of $W_k$ 
is called a chamber.
The moduli space 
$M_{B}(-\varphi_k(\alpha),\theta_k^{\dagger})$
does not change when $\theta_k^{\dagger}$
moves in a chamber. If 
two chambers $\Hat{C}_{\varphi_k(\alpha)}$ 
and $\Bar{C}_{\varphi_k(\alpha)}$ 
on $\varphi_k(\alpha)^{\perp}$ are adjacent
to each other having a common wall 
$w_k$, then for 
$\Hat{\theta}_k\in \Hat{C}_{\varphi_k(\alpha)}$, 
$\Bar{\theta}_k\in \Bar{C}_{\varphi_k(\alpha)}$ and 
$\Tilde{\theta}_k\in w_k$
we have a diagram:
\begin{equation}\label{thflip}
\xymatrix{
M_{B}(-\varphi_k(\alpha),
\Bar{\theta}_k)\ar[dr]^{f''}  & &
\ar@{-->}[ll]_{\kappa }
M_{B}(-\varphi_k(\alpha),
\Hat{\theta}_k). 
\ar[dl]_{f'} 
\\
& M_{B}(-\varphi_k(\alpha),
\Tilde{\theta}_k) 
& }
\end{equation}
Further, if both 
$M_{B}(-\varphi_k(\alpha),
\Hat{\theta}_k)$ and
$M_{B}(-\varphi_k(\alpha),
\Bar{\theta}_k)$ are
non-empty, then we see
that $f',f''$ are birational 
morphisms by general theory
of Thaddeus \cite{Th}. 

Within $\varphi_k(\alpha)^{\perp}$,
we have a chamber $C^{\PP^2}_{\varphi_k(\alpha)}$ such 
that $M_{\PP^2}(\ch(\alpha),H)\cong
M_{B}(-\varphi_k(\alpha),\theta_k)$ 
for any $\theta_k
\in C^{\PP^2}_{\varphi_k(\alpha)}$. 
In the case $\rk(\alpha)
=1$, $c_1(\alpha)=1$ and $\ch_2
(\alpha)=\frac{1}{2}-n$, 
diagrams~(\ref{thflip}) with $k=0,1$
give the two birational transformations
of the Hilbert schemes $(\PP^2)
^{[n]}$ (Theorem~\ref{wall-crossing}). 
In the case $\rk(\alpha)=r$,
$c_1(\alpha)=1$, $\ch_2(\alpha)=
\frac{1}{2}-n$ with arbitrary
$r>0$, we 
will describe these diagrams 
more explicitly in \cite{O}.

Similar phenomena as in 
(\ref{thflip}), sometimes called
Wall-crossing phenomena, 
occur by variation
of polarizations on
some surfaces $X$ in case of 
Gieseker-stability. 
However the polarization is 
essentially unique in our case $X=\PP^2$
since $\Pic\PP^2\cong\Z H$. So
our phenomena are of different
nature. We expect that Bridgeland
theory is useful to study such
phenomena systematically.
\vspace{0.7mm}\\

\noindent\emph{\bf{Convention}}\vspace{0.5cm}\\
Throughout this paper we work over $\C$. Any 
scheme is of finite type over $\C$. 
For a scheme $Y$, 
we denote by $\Coh(Y)$ the abelian category
of coherent sheaves on $Y$ and by $\D(Y)$ (respectively, 
$D^-(Y)$) the bounded (respectively, bounded above) 
derived category of $\Coh(Y)$.
For $E\in\Coh(Y)$, 
by $\dim E$ we denote the dimension of the 
support of $E$. For a ring $B$, by mod-$B$ 
we denote 
the abelian category of finitely generated 
right $B$-modules.
We denote by $\D(B)$ (respectively, $D^-(B)$)
the bounded (respectively, bounded above) derived 
category of mod-$B$.
For an abelian category $\A$ and a triangulated 
category $\T$, their Grothendieck groups are denoted 
by $K(\A)$ and $K(\T)$. For any object $E$ of $\A$
(resp. $\T$) we denote by $[E]$ the image of $E$ by 
the map $\A\to K(\A)$ (resp. $\T\to K(\T)$). 
When $\A=$ mod-$B$ and $\T=\D(Y)$, 
we simply write them $K(B)$ and $K(Y)$. 
For objects $E,F,G$ of $\T$, the
distinguished triangle 
$E\to F\to G\to E[1]$ is denoted by:
$$\xymatrix{
E \ar[rr]  & &F \ar[dl]\\
&  G\ar[ul]^{[1]} & 
}$$
For objects $F_0,\cdots,F_n$ in $\T$
we denote by $\langle
F_0,\cdots,F_n\rangle$ the
smallest full subcategory
of $\T$ containing $F_0,\cdots,F_n$,
which is closed under extensions.

\section{Generalities on 
Bridgeland stability conditions}
Here we collect some 
basic definitions and 
results of Bridgeland stability
conditions on triangulated categories
in \cite{B2}, \cite{B3}.

\subsection{Bridgeland 
stability conditions on 
triangulated 
categories}
Let $\A$ be an abelian  
category.
\begin{defn}\label{stf}
A stability function on $\A$ 
is a group homomorphism 
$Z\colon \K(\A)\to\C$ 
such that
$Z(E)\in\R_{>0}\,\exp
({\sqrt{-1}\pi\phi(E)})
\text{ with }0<\phi(E)
\leq 1$ for any nonzero
object $E$ of $\A$.
The real number $\phi(E)
\in(0,1]$ is called 
the phase of the object 
$E$.
A nonzero object $E$ of 
$\A$ is said to be
$Z$-(semi)stable if
for every proper subobject 
$0\neq F\subsetneq E$
we have $\phi(F)<\phi(E)$ 
(resp. $\le$).
\end{defn}

If we define the slope 
of $E$ by
$$\mu_{\sigma}(E) := 
-\frac{\mbox{Re}(Z(E))}
{\Im(Z(E))},$$
which possibly be infinity,
then a nonzero object $E$ of $\A$ 
is $Z$-(semi)stable if and only if 
$\mu_{\sigma}(F)<\mu_{\sigma}(E)$ 
($\text{resp.}~\le$) for any subobject 
$0\neq F\subsetneq E$ in $\A$.

The stability function $Z$ is said to have 
the Harder-Narasimhan property if
every nonzero object $E\in\A$ has
a finite filtration
\[0=E_0\subset E_1\subset \cdots\subset E_{n-1}
\subset E_n=E\]
whose factors $F_j=E_j/E_{j-1}$ are $Z$-semistable 
objects of $\A$ with
\[\phi(F_1)>\phi(F_2)>\cdots>\phi(F_n).\]

Let $\T$ be a triangulated category. We 
recall the definition of a t-structure
and its heart (cf. \cite{B2}).
\begin{defn}
A t-structure on $\T$ is a full subcategory
$\T^{\le 0}$ of $\T$ satisfying the following
properties.\\
$(1)$ $\T^{\le 0}[1]
\subset\T^{\le 0}$.\\ 
$(2)$ If one defines
$\T^{\ge 1}:=\left\{F\in\T\mid\Hom_{\T}
(G,F)=0 \text{ for any }G\in\T^{\le 0}
\right\},$ 
then for any object $E\in\T$ there is a 
distinguished triangle $G\to E\to F\to G[1]$
with $G\in\T^{\le 0}$ and $F\in\T^{\ge 1}$.
\end{defn}
We define $\T^{\le -i}:=\T^{\le 0}[i]$ and
$\T^{\ge -i}:=\T^{\ge 1}[i+1]$.
Then the heart of the t-structure is 
defined to be the full
subcategory $\A:=\T^{\le 0}\cap\T^{\ge 0}$.
It was proved in \cite{BBD} that 
$\A$ is an abelian category, with
the short exact sequences in $\A$
being precisely the triangles in $\T$
all of whose vertices are objects of
$\A$. A t-structure $\T^{\le 0}\subset\T$ is
said to be bounded if $$\T=\bigcup_{i,j\in\Z}
\T^{\le i}\cap\T^{\ge j}.$$
If $\A$ is the heart of a bounded t-structure
on $\T$, then we have $K(\A)\cong K(\T)$. 

\begin{defn}\label{stc}
A Bridgeland stability condition $\sigma$ 
on a triangulated category $\T$ is 
a pair $(Z,\A)$ of a group 
homomorphism $Z\colon K(\T)
\to\C$ and the heart $\A$ of 
a bounded t-structure on $\T$ 
such that $Z$ is a stability
function on $\A$ having the 
Harder-Narasimhan property.
\end{defn}

For each $n\in\Z$ and 
$\phi'\in(0,1]$, 
we define a full 
subcategory $\p(n+\phi')$ 
of $\T$ by
$$\p(n+\phi'):=\{E\in\T
\mid E[-n]\in\A \text{ is }
Z\text{-semistable} 
\text{ and }\phi(E[-n])
=\phi'\}.$$
For any $\phi\in\R$, 
a nonzero object 
$E$ of $\p(\phi)$ is said 
to be $\sigma$-semistable 
and $\phi$ is 
called the phase of $E$. 
$E\in\p(\phi)$ is said to be 
$\sigma$-stable if $\phi=n+\phi'$ 
with $n\in\Z$ and $\phi'
\in(0,1]$, and $E[-n]\in\A$ 
is $Z$-stable.
It is easy to see that each 
subcategory $\p(\phi)$ of $\T$ 
is an abelian category
(cf. \cite[Lemma~5.2]{B2}). 
$E\in\p(\phi)$ is 
$\sigma$-stable if and only if 
$E$ is a simple object 
in $\p(\phi)$. For any 
interval $I\subset\R$, $\p(I)$ 
is defined by 
$\p(I):=\langle\{\p(\phi)\mid 
\phi\in I\}\rangle.$ 
In particular the Harder-Narasimhan 
property implies that $\p((0,1])=\A$.
\begin{prop}\label{slice}
$(1)$ The pair $(Z,\p)$ of the
group homomorphism $Z\colon K(\T)\to\C$ and
the family $\p=\{\p(\phi)\mid\phi\in\R\}$
of full subcategories of $\T$ has
the following property.\\
$(\mathop{\mathrm{a}}\nolimits)$ 
$\p(\phi)$ is a full additive subcategory
of $\T$.\\
$(\mathop{\mathrm{b}}\nolimits)$ 
$\p(\phi+1)=\p(\phi)[1]$.\\
$(\mathop{\mathrm{c}}\nolimits)$ 
If $\phi_1>\phi_2$ and $E_i\in\p(\phi_i)$,
then $\Hom_{\T}(E_1,E_2)=0$.\\
$(\mathop{\mathrm{d}}\nolimits)$ 
$Z(E)\in\R_{>0}\exp(\sqrt{-1}\pi\phi)$
for any nonzero object $E$ of $\p(\phi)$.\\
$(\mathop{\mathrm{e}}\nolimits)$ 
For a nonzero object $E\in\T$, we have
a collection of triangles
$$\xymatrix{
0=E_0 \ar[rr]  & &E_1 
\ar[dl] \ar[rr] & & E_2 
\ar[r]\ar[dl] 
& \cdots \ar[rr] & & E_n
 =E \ar[dl]\\
&  F_1 \ar[ul]^{[1]} & &
 F_2 \ar[ul]^{[1]}& & & 
F_n \ar[ul]^{[1]}&
}$$
such that $F_j\in\p(\phi_j)$ with 
$\phi_1>\phi_2>\cdots>\phi_n$.\\
$(2)$ Giving a stability condition 
$\sigma=(Z,\A)$ on $\T$ is
equivalent to giving a pair $(Z,\p)$ 
with the above properties.
\end{prop}
\proof
See 
\cite[Definition~5.1 and 
Proposition~5.3]{B2}. 
Originally the pair $(Z,\p)$
is called the stability 
condition
$\sigma$ in \cite{B2}.
\endproof

The filtration in (e) of Proposition~\ref{slice}
is called the Harder-Narasimhan filtration
of $E$ and the objects $F_j$ are called 
$\sigma$-semistable factors of $E$.
We can easily check that the Harder-Narasimhan
filtration is unique up to isomorphism.
For a Bridgeland stability condition
$\sigma=(Z,\A)$ (or $(Z,\p)$), $Z$,
$\A$ and $\p$ is denoted by 
$Z_\sigma$, $\A_\sigma$ and $\p_\sigma$.

\subsection{Bridgeland stability conditions 
on smooth projective surfaces}
\label{brionsurf}
Let $X$ be a smooth complex projective surface. 
The Chern character of an object $E$ of 
$\D(X)$ is the element 
$\ch(E):=(\rk(E),c_1(E),\ch_2(E))$
of the lattice
$\N(X):=\Z\oplus\ns(X)\oplus\frac{1}{2}\Z$.
We define the Euler form on the Grothendieck 
group $K(X)$ of $X$ by
\begin{equation}
\chi(E,F):=\Sigma_i (-1)^i\dim_{\C}\Hom_{\D(X)}(E,F[i]).
\end{equation}
Let $K(X)^{\perp}=\{\alpha\in K(X)\mid\chi(\alpha,\beta)=0 
\ \text{for each} \ \beta\in\K(X)\}$
and $K(X)/K(X)^{\perp}$ is called 
the \emph{numerical Grothendieck group} of $\D(X)$.\\
By the Riemann-Roch theorem the Chern character 
gives an inclusion 
$K(X)/K(X)^{\perp}\to \N(X)$. 
Furthermore we define a symmetric bilinear 
form $(\ \cdot\ ,\ \cdot\ )_M$ on $\N(X)$, called Mukai
pairing, by the following formula
\begin{equation}\label{inp}
((r_1,D_1,s_1),(r_2,D_2,s_2))_M
:=D_1\cdot D_2-r_1s_2-r_2s_1.
\end{equation}
This bilinear form makes $\N(X)$ 
a lattice of signature $(2,\rho)$
by the Hodge Index Theorem, 
where $\rho\ge1$ is the 
Picard number of $X$. 

A Bridgeland stability condition $\sigma=(Z,\A)$ 
is said to be \emph{numerical} if 
there is a vector $\pi(\sigma)\in\N(X)\otimes\C$
such that
\begin{equation}\label{pi}
Z(E)=(\pi(\sigma),\ch(E))_M
\end{equation}
for any $[E]\in K(X)$. $\sigma$ is
said to be \emph{local finite} if
it satisfies some technical 
conditions~\cite[Definition 5.7]{B2}.

The set of all the numerical local finite
Bridgeland stability conditions 
on $\D(X)$ is denoted by $\Stab (X)$.
It is shown in~\cite[Section 6]{B2} that 
$\Stab (X)$ has a natural structure as a 
complex manifold. The map 
\begin{equation}\label{pim}
\pi\colon\Stab(X)\to\N(X)\otimes\C,
\end{equation}
defined by (\ref{pi}), is holomorphic.

For the fixed heart $\A$ of a bounded 
t-structure on $\D(X)$, we write 
$$\Stab(\A):=\{\sigma\in\Stab(X)
\mid \A_{\sigma}=\A\}.$$

\subsection{$\grp$ action 
on $\Stab(X)$}
\label{grp}
Let $\widetilde{\GL}^{+}
(2, \mathbb{R})$ be the 
universal cover of $\GL ^{+}(2, \mathbb{R})
=\{T\in\GL(2,\R)\mid\det T>0\}$.  
The group $\grp$ can be viewed as the 
set of pairs $(T,f)$
where $T\in\GL^+(2,\R)$ and $f$ 
is the automorphism of $\R\cong\widetilde{
S^1}$ such that $f$ covers the automorphism
$\Bar{T}$ of $S^1\cong\left(\R^2\setminus{0}
\right)/\R_{>0}$ induced by $T$.

The topological space $\Stab (X)$ carries 
the right action of 
the group $\widetilde{\GL}^{+}(2, \mathbb{R})$ 
\cite[Lemma 8.2]{B2} as follows. 
Given $\sigma\in
\Stab(X)$ and $g=(T,f)\in\grp$, 
a new stability 
condition $\sigma g$ is defined to be 
the pair $(Z_{\sigma g},\p_{\sigma g})$ 
where 
$Z_{\sigma g}:=T^{-1}\circ Z_{\sigma}$ and 
$\p_{\sigma g}(\phi):=\p_{\sigma}(f(\phi))$
for $\phi\in\R$, 
where we identify $\C$ 
with $\R^2$ by 
$$x+\sqrt{-1}y\mapsto 
\begin{pmatrix}
x\\
y
\end{pmatrix}.$$
It is easy to check that the pair 
$(Z_{\sigma g}, \p_{\sigma g})$ satisfies 
the properties
of Proposition~\ref{slice}~(1). Hence by
Proposition~\ref{slice}~(2), we have
$\sigma g=(Z_{\sigma g},\p_{\sigma g})
\in\Stab(X)$.
We remark that the sets of 
the (semi)stable objects of 
$\sigma$ and $\sigma g$ are the same, 
but the phases have 
been relabelled. For our purpose,
it is convenient to 
introduce the following definition.
\begin{defn}
Two stability conditions $\sigma, \sigma'
\in\Stab(X)$ are said to be 
$\grp$-equivalent to each other if 
$\sigma$ and $\sigma'$ are in a single
$\grp$ orbit. 
\end{defn}

For any element $T\in\GL^+(2,\R)$,
the right $\GL^+(2,\R)$ action on 
$\N(X)\otimes\C$ is defined by 
$\id_{\N(X)}\otimes T^{-1}$.
Hence the $\grp$ acts on $\N(X)\otimes\C$ via 
the covering map 
$$\grp\to\GL^+(2,\R)\colon
(T,f)\mapsto T.$$ The map
$\pi\colon\Stab(X)\to\N(X)\otimes\C$
is equivariant for these $\grp$ actions.

\subsection{Moduli functors of 
Bridgeland semistable objects}
\label{modu}
For $\sigma=(Z,\A)\in\Stab(X)$ 
and $\alpha\in K(X)$, we define 
a moduli functor 
$$\M_{\D(X)}(\ch(\alpha),\sigma)\colon
(\text{scheme}/\C)\to(\text{sets})
\colon S\mapsto \M_{\D(X)}(\ch(\alpha),
\sigma)(S)$$ as follows, where 
$(\text{scheme}/\C)$ is the category 
of schemes of finite type over $\C$
and $(\text{sets})$ is the category
of sets. For a 
scheme $S$, the set 
$\M_{\D(X)}(\ch(\alpha),\sigma)(S)$ 
consists of isomorphism 
classes of $E\in \D(X\times S)$ 
such that for every closed point 
$s\in S$ the restriction to the 
fiber
\begin{equation*}
E_s:=\mathbf L \iota_{X
\times \{s\}}^{\ast}E
\end{equation*}
is a $\sigma$-semistable object 
in $\A$ with $\ch(E_s)=\ch
(\alpha)\in \N(X)$, 
where $\iota_{X\times \{s\}}$ 
is the embedding 
$$\iota_{X\times \{s\}}\colon 
X\times\big\{s\big\}
\to X\times S.$$

\noindent Note that by 
definition each object $E_s$ 
belongs to $\A \subset
\D(X)$ for every closed point 
$s\in S$, so $\ch(E_s)\in\N(X)$
is well-defined. Let $\M^s_{\D(X)}
(\ch(\alpha),\sigma)$ be the 
subfunctor of $\M_{\D(X)}(\ch(\alpha),
\sigma)$ corresponding to 
$\sigma$-stable objects of $\A$.

Since the action of $\grp$ does 
not change the set of 
(semi)stable objects, for any 
$g\in\grp$ there exists an integer 
$n$ such that the shift 
functor $[n]$ gives an isomorphism
\begin{equation}\label{gisom}
\M_{\D(X)}(\ch(\alpha),\sigma)\cong
\M_{\D(X)}((-1)^n\ch(\alpha),\sigma g)
\colon E\mapsto E[n].
\end{equation}

Here we recall the definition 
of a moduli space. For a scheme 
$Z$, we denote by $\underline{Z}$ 
the functor 
$$\underline{Z}\colon 
(\text{scheme}/\C)\to(\text{sets})
\colon S\mapsto \Hom(S,Z).$$
The Yoneda lemma tells us that 
every natural transformation 
$\underline{Y}\to\underline{Z}$ 
is of the form $\underline{f}$ 
for some morphism $f\colon Y\to Z$ 
of schemes, where $\underline{f}$ 
sends $t\in\underline{Y}(T)$ to 
$\underline{f}(t)=f\circ t\in
\underline{Z}(T)$ for any scheme $T$.
A functor $(\text{scheme}/\C)
\to(\text{sets})$ isomorphic to 
$\underline{Z}$ is said to be 
represented by $Z$.

In the terminology introduced 
by Simpson \cite[Section~1]{S},
a \emph{moduli space} is a scheme 
which 'corepresents' 
a moduli functor. 

\begin{defn}\label{corep}
Let $\M\colon(\text{scheme}/
\C)\to(\text{sets})$
be a functor, $M$ a scheme and 
$\psi\colon\M\to \underline{M}$ 
a natural transformation. 
We say that $M$ corepresents $\M$ 
if for each scheme $Y$ and each 
natural transformation 
$h\colon\M\to\underline{Y}$, there 
exists a unique morphism $\sigma
\colon M\to Y$ such that 
$h=\underline{\sigma}\circ\psi$:
$$\xymatrix{
\M \ar[dr]^{h}\ar[d]_{\psi}&\\
\underline{M}\ar[r]_{\underline
{\sigma}}&\underline{Y}
}$$
\end{defn}
This characterizes $M$ up to a 
unique isomorphism. If $M$ 
represents $\M$ we say that $M$ 
is a fine moduli space.

For any functor $\M\colon
(\text{scheme}/\C)\to(\text{sets})$,
we consider the sheafication 
of $\M$
$${}^{sh}\M\colon(\text{scheme}
/\C)\to(\text{sets})$$ 
with respect to the
Zariski topology. 
For a scheme $S$, ${}^{sh}
\M(S)$ is defined as follows.
For an open cover 
$\mathcal{U}=\{
U_i\}$ of $S$, $S=\cup U_i$, let $\M_{
\mathcal{U}}:=\{(E_i)\in\prod\M
(U_i)\mid E_i|_{U_i\cap U_j}=
E_j|_{U_i\cap U_j}\}$.
If $\mathcal{V}$ is a refinement of 
$\mathcal{U}$, then we have a 
natural map $\M_{\mathcal{U}}\to
\M_{\mathcal{V}}$. The set of 
open covers forms a direct system
with respect to the preorder 
defined by refinement.
${}^{sh}\M(S)$ is defined by
\begin{equation}
\label{shi}
{}^{sh}\M(S):=\lim_{\substack{
\xrightarrow{\hspace{
0.3cm}}\\
\mathcal{U}}}\M_{\mathcal{U}}.
\end{equation}
Actually, the limit can be 
computed over affine coverings
only, because every covering 
$\mathcal{U}$ has a refinement
which is affine.
Since any scheme $Y$ satisfies
$\underline{Y}\cong{}^{sh}
\underline{Y}$, we have 
\begin{equation}\label{sh}
\Hom(\M,\underline{Y})\cong\Hom(
{}^{sh}\M,\underline{Y}).
\end{equation}
In particular, a scheme $M$ 
corepresents $\M$ if and only if 
$M$ corepresents ${}^{sh}\M$.

\section{Geometric Bridgeland 
stability conditions}
Let $X$ be a smooth projective 
surface. In this section, we 
introduce the notion of geometric 
Bridgeland stability conditions 
on $\D(X)$ and see that if 
$\sigma\in\Stab(X)$ is 
geometric, then under suitable 
assumptions the above functor 
$\M^s_{\D(X)}(\ch(\alpha),\sigma)$  
$\left(\text{resp.}~\M_{\D(X)}
(\ch(\alpha),\sigma)~\right)$ is 
corepresented by a certain moduli 
space of Gieseker-(semi)stable 
coherent sheaves on 
$X$.

\subsection{Twisted 
Gieseker-stability and $\mu$-stability}
We recall the notion of twisted 
Gieseker-stability and 
$\mu$-stability. For details, 
we can consult \cite{HL}, 
\cite{MW}.
Take $\gamma,\omega\in\ns(X)
\otimes\R$, 
and suppose that $\omega$ 
is in the ample cone
$$\Amp(X)=\big\{\omega\in\ns(X)
\otimes\R\mid \omega^2>0 
\text{ and }\omega\cdot C>0
\text{ for any curve }
C\subset X\big\}.$$
For a coherent sheaf $E$ with 
$\rk(E)\neq0$, 
define  
$\mu_{\omega}(E)$ and 
$\nu_{\gamma}(E)$ by
\begin{equation}\label{hil1}
\mu_{\omega}(E):=\frac{c_1(E)\cdot\omega}
{\rk(E)},\hspace{1cm} 
\nu_{\gamma}(E):=\frac{\ch_2(E)}{\rk(E)}-
\frac{c_1(E)\cdot K_X}{2\rk(E)}
-\frac{c_1(E)\cdot\gamma}{\rk(E)}.
\end{equation}

\begin{defn}Let 
$E$ be a torsion free sheaf.\\
$($\emph{i}$)$ $E$ is said to be 
$(\gamma,\omega)$-semistable if 
for every proper nonzero subsheaf $F$
of $E$ we have 
\begin{equation}\label{hil2}
\left(\mu_{\omega}(F), \nu_{\gamma}(F)\right)\le
\left(\mu_{\omega}(E), \nu_{\gamma}(E)\right)
\end{equation}
in the lexicographic order, 
namely $\mu_{\omega}(F)<\mu_{\omega}(E)$
or $\mu_{\omega}(F)=\mu_{\omega}(E),
\nu_{\gamma}(F)\le\nu_{\gamma}(E)$.
$E$ is said to be $(\gamma,\omega)$-stable
if $\left(\mu_{\omega}(F), \nu_{\gamma}(F)\right)<
\left(\mu_{\omega}(E), \nu_{\gamma}(E)\right)$
for any such $F$.\\
$($\emph{ii}$)$ $E$ is said to be $\mu_{\omega}$-semistable 
if $\mu_{\omega}(F)\le\mu_{\omega}(E)$
for any such $F$. 
$E$ is said to be 
$\mu_{\omega}$-stable if in addition 
$\mu_{\omega}(F)<\mu_{\omega}(E)$
for any $F$ with $\rk F<\rk E$.
\end{defn}  

$(\gamma,\omega)$-stability is
called twisted Gieseker-stability
in \cite{To}.
Correspondingly to these semistability
notions, every torsion free sheaf 
$E$ on $X$ has a unique Harder-Narasimhan 
filtration (cf. \cite[Example~4.16 and 4.17]{J}). If
$$0=E_0\subset E_1\subset \cdots 
\subset E_{n-1}\subset E_n=E$$
is the Harder-Narasimhan 
filtration with respect to 
$\mu_{\omega}$-semistability, 
we define $\mu_{\omega\text{-min}}(E)
:=\mu_{\omega}(E_n/E_{n-1})$ and 
$\mu_{\omega\text{-max}}(E):=
\mu_{\omega}(E_1)$. 

\begin{thm}\label{bgine}
\emph{\bf{(Bogomolov-Gieseker 
Inequality)}}. 
Let $X$ be a smooth projective 
surface and $\omega$ an ample 
divisor on $X$. If $E$ is a 
$\mu_{\omega}$-semistable 
torsion free sheaf on 
$X$, then
\begin{equation*}\label{BG}
c_1^2(E)-2\rk(E)\ch_2(E)\ge 0.
\end{equation*}
\end{thm}
\proof
See \cite[Theorem~3.4.1]{HL}.
\endproof

We take $\alpha\in K(X)$ with 
$\rk(\alpha)>0$ and consider the moduli functor 
$\M_{X}(\ch(\alpha),\gamma,\omega)$
of $(\gamma,\omega)$-semistable 
torsion free sheaves $E$ 
with $\ch(E)=\ch(\alpha)\in\ns(X)$.
Let $\M^s_{X}(\ch(\alpha),\gamma,\omega)$ 
be the subfunctor of $\M_{X}(\ch(\alpha),
\gamma,\omega)$
corresponding to 
$(\gamma,\omega)$-stable 
ones.

We denote by $M_{X}(\ch(\alpha),\gamma,\omega)$
the moduli space of 
$(\gamma,\omega)$-semistable 
torsion-free sheaves if it exists.
When $\omega$ is an integral ample
divisor and $\gamma\in\ns(X)\otimes\Q$,
the moduli space $M_{X}(\ch(\alpha),
\gamma,\omega)$ exists \cite[Theorem~5.7]{MW}. 
Furthermore if $\gamma=0$, we write $M_{X}(\ch(\alpha),\omega)$
instead of $M_{X}(\ch(\alpha),0,\omega)$ for
the sake of simplicity. In this case 
there is an open subset $M^s_{X}(\ch(\alpha),\omega)$ of 
$M_{X}(\ch(\alpha),\omega)$ that corepresents 
the functor $\M^s_{X}(\ch(\alpha),\omega)$ 
\cite[Theorem~4.3.4]{HL}.

\subsection{Geometric 
Bridgeland stability conditions}
\label{ABL}

We construct some Bridgeland 
stability conditions on $\D(X)$ 
following \cite{ABL}.
For every coherent sheaf $E$ on 
$X$, we denote the torsion part 
of $E$ by $E_{\tor}$ and the 
torsion free part of $E$ by 
$E_{\fr}=E/E_{\tor}$. Suppose 
that $\beta,\omega \in \ns(X)
\otimes\R$ with $\omega\in
\Amp(X)$, then we define two 
full subcategories $\mathfrak{T}$ 
and $\mathfrak{F}$ of $\Coh(X)$ 
as follows;
\begin{equation*}
\begin{split}
\text{ob}(\mathfrak{T})&=\{
\text{torsion sheaves}\}\cup
\{E\mid E_{\fr}\neq 0 \text{ and } 
\mu_{\omega\text{-min}}(E_{\fr})>
\beta\cdot\omega\}\\
\text{ob}(\mathfrak{F})&=\{E\mid 
E_{\tor}=0\text{ and } \mu_{
\omega\text{-max}}(E)\le \beta
\cdot\omega\}.
\end{split}
\end{equation*} 
We define a pair $\sigma_{(\beta,\omega)} = 
(Z_{(\beta,\omega)},\A_{(\beta,\omega)})$ 
of the heart $\A_{(\beta,\omega)}$ of 
a bounded t-structure on $\D(X)$ 
and a stability function $Z_{(\beta,\omega)}$
on $\A_{(\beta,\omega)}$ in the following way.\\

\begin{defn}\label{def of geometric}
A full subcategory $\A_{(\beta,\omega)}$ 
of $\D(X)$ is defined as follows;
$$\A_{(\beta,\omega)}:=\{E\in\D(X)\mid\
\mathcal{H}^i(E)=0 \text{ for all } i \neq 0,1 
\text{ and } \mathcal{H}^0(E)\in\mathfrak{T}
\text{ and } \mathcal{H}^{-1}(E)\in\mathfrak{F}\}.$$
The group homomorphism $Z_{(\beta,\omega)}$ 
is defined by $Z_{(\beta,\omega)}(\alpha):=
(\exp(\beta+\sqrt{-1}\omega),\ch(\alpha))_M$, 
where $$\exp(\beta+\sqrt{-1}\omega)=
\left(1,\beta+\sqrt{-1}\omega,\frac{1}{2}(\beta^2-\omega^2)
+\sqrt{-1}(\beta\cdot\omega)\right) \in \N(X)$$
and $(\ \cdot\ ,\ \cdot\ )_M$ is the Mukai pairing defined in
\S~$\ref{brionsurf}$.
\end{defn}
From the general theory called tilting
we see that $\A_{(\beta,\omega)}$ is the 
heart of a bounded t-structure on $\D(X)$
(for example, see \cite[\S~3]{B2}).
By definition, for $\alpha\in K(X)$ with 
$\ch(\alpha)=(r,c_1,\ch_2)$ we have
\begin{equation}\label{can}
Z_{(\beta,\omega)}(\alpha)=-\ch_2+c_1\cdot\beta+\frac{r}{2}
(\omega^2-\beta^2)+\sqrt{-1}\omega\cdot(c_1-r\beta). 
\end{equation}
Furthermore if $r\neq 0$, we can write 
\begin{equation}\label{cha}
Z_{(\beta,\omega)}(\alpha)=\frac{1}{2r}((c_1^2-2r \ch_2)+
r^2 \omega^2-(c_1-r\beta)^2)+\sqrt{-1}\omega(c_1-r\beta).
\end{equation}

Our $\sigma_{(\beta,\omega)}$
is slightly different from 
that in \cite{B3}, \cite{To}.

\begin{prop}\emph{\bf{\cite[Corollary 2.1]{ABL}}}
\label{abl}For each pair 
$\beta,\omega \in \ns(X)\otimes\Q$ with 
$\omega\in\Amp(X)$, $\sigma_{(\beta,\omega)}$ is a 
Bridgeland stability condition on $\D(X)$.
\end{prop}

For general $\beta,\omega\in
\ns(X)\otimes\R$, we do not
know whether $\sigma_{(\beta,
\omega)}$ belongs to $\Stab(X)$
or not since we do not know if 
$Z_{(\beta,\omega)}$ has the 
Harder-Narasimhan property. 
If $\beta,\omega\in\ns(X)
\otimes\Q$ it directly follows 
from
\cite[Proposition~7.1]{B3}.
However we consider the 
following definition.
\begin{defn}\label{deg}
$\sigma\in\Stab(X)$ is called 
geometric if $\sigma$ is 
$\grp$-equivalent to $\sigma_
{(\beta,\omega)}$ for some 
$\beta,\omega\in\ns(X)\otimes
\R$ with $\omega\in\Amp(X)$.
\end{defn}  

We have the following criterion 
due to \cite{B3} 
for $\sigma\in\Stab(X)$
to be geometric. It reduces the
proof of Theorem~\ref{our}
to easy calculations 
(\S~\ref{proof}).
\begin{prop}\label{geometric1}
$\sigma \in \Stab(X)$ is 
geometric if and 
only if \\
$1$. For all $x \in X$, 
the structure sheaves $\mo_x$ 
are $\sigma$-stable of the same phase.\\
$2$. There exist $T\in\GL^+(2,\R)$ and
$\beta,\omega\in\ns(X)\otimes\R$ such that 
$\omega^2>0$ and $$\pi(\sigma)T
=\exp(\beta+\sqrt{-1}\omega),$$
where $\pi\colon\Stab(X)\to\N(X)$ is
defined by $(\ref{pim})$ and 
$\GL^+(2,\R)$ action on $\N(X)\otimes\C$
is defined in \S~$\ref{grp}$.
\end{prop}
\proof
From \cite[Lemma~10.1 and Proposition~10.3]{B3} 
the assertion holds because \cite[Lemma~6.3
and Lemma~10.1]{B3}
hold for an arbitrary smooth projective surface. 
However we give the proof of this proposition
for the reader's convenience.


The only if part is easy.  
By \cite[Lemma~6.3]{B3}, 
for any closed point 
$x\in X$ the structure sheaf $\mo_x$ 
is a simple object of the abelian category 
$\A_{(\beta,\omega)}$, hence 
$\sigma_{(\beta,\omega)}$-stable for any 
$\beta,\omega\in\ns(X)$ with $\omega\in\Amp(X)$. 
Since $\grp$ action does not change 
stable objects, $\mo_x$ is also $\sigma$-stable.
Furthermore since the map $\pi$ is 
equivariant for $\grp$ actions,
$\sigma$ also satisfies condition~2 
(cf. \S~\ref{grp}). 

Now we consider the if part. 
We show that $\sigma g=\sigma_{
(\beta,\omega)}$ for some 
$g=(T,f)\in\grp$, where $\beta,
\omega$ and $T$ are as in the
condition~2. We may assume 
$\pi(\sigma)=
\exp(\beta+\sqrt{-1}\omega)$ for some $\beta,\omega
\in\ns(X)\otimes\R$ with $\omega^2>0$. The kernel of 
the homomorphism $\grp\to\GL^+(2,\R)$ acts 
on $\Stab(X)$ by even shifts, so we may assume 
furthermore that $\mo_x\in\p_\sigma(1)$ for all 
$x\in X$. 

We show that $\omega$ is ample. 
It is enough to show that 
$C\cdot\omega>0$ for any curve $C\subset X$. 
The condition~1 and 
\cite[Lemma~10.1(c)]{B3} 
show that the torsion sheaf $\mo_C$ 
lies in the subcategory $\p_\sigma((0,1])$. 
If $Z_\sigma(\mo_C)$ lies on 
the real axis it follows that 
$\mo_C\in\p_\sigma(1)$ which is impossible 
by \cite[Lemma~10.1(b)]{B3}. 
Thus $\Im Z_\sigma(\mo_C)=C\cdot\omega>0$. 

The same argument of STEP 2 in 
\cite[Proposition~10.3]{B3} 
holds and we see that 
$\p_\sigma((0,1])=\A_{(\beta,\omega)}$.
\endproof

\subsection{Moduli spaces 
corepresenting 
$\M_{\D(X)}(\ch(\alpha),
\sigma_{(\beta,\omega)})$ and 
$\M^s_{\D(X)}(\ch(\alpha),
\sigma_{(\beta,\omega)})$}
\label{moco}
In this subsection we fix 
$\alpha\in K(X)$ with 
$\ch(\alpha)=(r,c_1,\ch_2)
\in\N(X)$, $r>0$ and $\beta 
\in\ns(X)\otimes\mathbb{R}$, 
$\omega\in\ns(X)$ 
with $\omega$ ample. We put 
\begin{equation}\label{ep}
\e:=\Im Z_{(\beta,\omega)}
(\alpha)=c_1\cdot\omega-r
\beta\cdot\omega\in\R
\end{equation}
and $\gamma:=\beta-
\frac{1}{2}K_X\in\ns(X)
\otimes\mathbb{R}$. 
We take $0<t\le1$ and 
assume that $\sigma_{(\beta,
t\omega)}=(Z_{(\beta,t\omega)},
\A_{(\beta,t\omega)})$ satisfies 
the Harder-Narasimhan 
property, that is, $\sigma_{
(\beta,t\omega)}\in\Stab(X)$.
We will show that if $\e>0$ 
is small enough and the moduli space 
$M_{X}(\ch(\alpha),\gamma,\omega)$ 
exists, then it corepresents the moduli functor
$\M_{\D(X)}(\ch(\alpha),\sigma_{
(\beta,t\omega)}).$

\begin{lem}\label{qlim}
For any $\sigma_{(\beta,t\omega)}$-semistable 
object $E \in \A_{(\beta,t\omega)}$ 
with $[E]=\alpha$, the following hold. 
\\
$(1)$ Assume that $0<\e\le t$ and 
$\Re Z_{(\beta,t\omega)}(\alpha)\ge 0$. 
Then $E$ is a torsion free sheaf.\\
$(2)$ Furthermore assume that $\e\le\frac{1}{r}$.  
Then $E$ is a $\mu_{\omega}$-semistable torsion free sheaf.
\end{lem}
\proof
$(1)$ For a contradiction we assume 
that $\mathcal{H}^{-1}(E)\neq 0$ 
and take $\ch(\mathcal{H}^{-1}(E))
=(r',c_1',\ch_2')\in\N(X)$. 
Then there exists an exact sequence 
in $\A_{(\beta,t\omega)}$, 
\begin{equation}\label{ex1}
0 \to \mathcal{H}^{-1}(E)[1]\to 
E\to \mathcal{H}^{0}(E)\to 0
\end{equation}
and we have $Z_{(\beta,t\omega)}
(E)=Z_{(\beta,t\omega)}
(\mathcal{H}^0(E))+Z_{(\beta,
t\omega)}(\mathcal{H}^{-1}(E)
[1])$. 
Since $\Im Z_{(\beta,t\omega)}
(\mathcal{H}^0(E))>0$ and
$\Im Z_{(\beta,t\omega)}
(\mathcal{H}^{-1}(E)[1])\ge0$,
we get 
$$0\le t\omega \cdot (-c_1'+
r'\beta)=\Im Z_{(\beta,t\omega)}
(\mathcal{H}^{-1}(E)[1]) <
\Im Z_{(\beta,t\omega)}(E)=t\e.$$
By the Hodge Index Theorem, we have
\begin{equation}\label{ineq}
(-c_1'+r'\beta)^2<\frac{\e^2}
{\omega^2}\le t^2.
\end{equation} 

Here we assume that $\mathcal{H}^{-1}
(E)$ is $\mu_{\omega}$-semistable. 
Then by Theorem~\ref{bgine} 
we have $-({c'_1}^2-2r'\ch_2')\le 0$. 
It follows from (\ref{cha}), (\ref{ineq}) 
and ${r'}^2\omega^2\in\Z_{>0}$ that
\begin{equation*}
  \begin{split}
\Re Z_{(\beta,t\omega)}(\mathcal{H}^{-1}(E)[1])&=
\frac{1}{2r'}(-(c_1'^2-2r'\ch_2')-r'^2t^2\omega^2+
(c_1'-r'\beta)^2)\\
&<\frac{1}{2r'}(-r'^2\omega^2+1)t^2\le 0.
\end{split}
\end{equation*}
In the general case, $\mathcal{H}^{-1}(E)$ factors into 
$\mu_{\omega}$-semistable sheaves and 
we also get the inequality 
$$\Re Z_{(\beta,t\omega)}(\mathcal{H}^{-1}(E)[1])<0.$$ 
Hence we have $0<\mu_{\sigma_{(\beta,t\omega)}}
(\mathcal{H}^{-1}(E)[1])$.

On the other hand by the assumption that 
$\Re Z_{(\beta,t\omega)}(E)\ge0$, we have 
$\mu_{\sigma_{(\beta,t\omega)}}(E)\le 0$.
Thus we have $\mu_{\sigma_{(\beta,t\omega)}}(E)
<\mu_{\sigma_{(\beta,t\omega)}}(\mathcal{H}^{-1}(E)[1])$.
This contradicts the fact that 
$E$ is $\sigma_{(\beta,t\omega)}$-semistable since 
$\mathcal{H}^{-1}(E)[1]$ is a subobject of $E$ in 
$\A_{(\beta,t\omega)}$ by (\ref{ex1}). 
Thus $\mathcal{H}^{-1}(E)=0$ and $E$ is
a sheaf. 

Next we show that $E$ is torsion free. 
We assume that 
$E$ has a torsion $E_{\tor}\neq 0$. 
In the case $\dim E_{\tor} =1$, 
we have $m:=\omega \cdot 
c_1(E_{\tor})\ge 1$.
Since $E\in\A_{(\beta,t\omega)}$ we get 
$t\omega\cdot \beta < 
\mu_{t\omega}(E_{\fr})=\frac{tc_1\cdot\omega-mt}{r}$. 
However by (\ref{ep}),
$t\omega\cdot\beta=\frac{tc_1\cdot\omega-t\e}{r}$.
This implies that $\e>m\ge 1$. This contradicts 
the assumption that $\e\le t\le 1$. 
In the case $\dim E_{\tor} =0$, we get a nonzero 
subobject $E_{\tor}$ of $E$ in $\A_{(\beta,t\omega)}$.
However the slope $\mu_{\sigma_{(\beta,t\omega)}}
(E_{\tor})$ is infinity and greater than 
$\mu_{\sigma_{(\beta,t\omega)}}(E)$. 
This contradicts the fact that $E$ is 
$\sigma_{(\beta,t\omega)}$-semistable. \\
$(2)$ By $(1)$, $E$ is a torsion free sheaf. 
For a contradiction we assume that 
$E$ is not $\mu_{\omega}$-semistable. 
Then there exists an exact 
sequence in $\Coh(X)$ 
\begin{equation*}
0\to E''\to E \to E'\to 0.
\end{equation*} 
Here $E'$ is a $\mu_{\omega}$-semistable 
factor of $E$ 
with the smallest slope $\mu_{\omega}(E')$.
Since $E\in\A_{(\beta,t\omega)}$, we have 
$t\omega\cdot \beta < \mu_{t\omega\text{-min}}(E)=
\mu_{t\omega}(E')$.
Hence
\begin{equation*}
\mu_{t\omega}(E)-\mu_{t\omega}(E')<\mu_{t\omega}(E)-
t\omega\cdot \beta=t\e/r.
\end{equation*}
On the other hand, since 
$\mu_{\omega}(E)-\mu_{\omega}(E')>0$
and $\rk(E')c_1\cdot\omega-rc_1(E')\cdot
\omega$ is an integer, we have
\begin{equation*}
\mu_{\omega}(E)-\mu_{\omega}(E') =
\frac{\rk(E')c_1\cdot\omega-rc_1(E')\cdot
\omega}{r\rk(E')}> 1/r^2.
\end{equation*}
Hence we get $\e/r >\mu_{\omega}(E)-\mu_{\omega}(E')
> 1/r^2$ and 
this contradicts the assumption that $\e\le\frac{1}{r}$. 
Thus $E$ is $\mu_{\omega}$-semistable.
\endproof

Next we consider the relationship 
between $\sigma_{(\beta,t\omega)}$ 
and the $(\gamma,\omega)$-stability,
where $\gamma=\beta-\frac{1}{2}K_X$. 
By (\ref{can}) the slope 
$\mu_{\sigma_{(\beta,t\omega)}}(E)$ 
is written as 
\begin{equation}\label{sig1}
\mu_{\sigma_{(\beta,t\omega)}}(E)=
\frac{\nu_{\gamma}(E)-
\frac{1}{2}(t^2\omega^2-\beta^2)}
{t\mu_{\omega}(E)-t\beta\cdot\omega}
\end{equation}
for any coherent sheaf $E\in\Coh(X)$ 
with $\rk(E)\neq0$.

\begin{thm}
\label{limit}
Assume that $0<\e\le\min\{t,
\frac{1}{r}\}$ and $\Re Z_{(
\beta,t\omega)}(\alpha)\ge 0$.
Then for $E \in \A_{(\beta,
t\omega)}$ with $[E]=\alpha$, 
$E$ is $\sigma_{(\beta,
t\omega)}$-$($semi$)$stable 
if and only if $E$ is a $(\gamma,
\omega)$-$($semi$)$stable 
torsion free sheaf.
\end{thm}
\proof
$\Rightarrow$) From Lemma~\ref{qlim}, $E$ is 
a $\mu_{\omega}$-semistable 
torsion free sheaf. Hence to see that $E$ is 
$(\gamma,\omega)$-(semi)stable 
it is enough to show that 
for any subsheaf $F\subset E$ 
with $E/F$ torsion free and 
$\mu_{\omega}(F)=\mu_{\omega}(E)$,  
the inequality $\nu_{\gamma}(F)<
\nu_{\gamma}(E), \ (\text{resp. }\le)$ holds. 
Since $E$ is $\mu_{\omega}$-semistable 
and $\mu_{\omega}(F)=
\mu_{\omega}(E/F)=\mu_{\omega}(E)$, 
both $F$ and $E/F$ are 
$\mu_{\omega}$-semistable and belong 
to $\A_{(\beta,t\omega)}$. 
Hence the exact sequence in $\Coh(X)$ 
$$0\to F\to E\to E/F\to 0$$
is also exact in $\A_{(\beta,t\omega)}$. 

Since $E$ is $\sigma_{(\beta,t\omega)}$-(semi)stable, 
we have $\mu_{\sigma_{(\beta,t\omega)}}(F)<
\mu_{\sigma_{(\beta,t\omega)}}(E)$, 
$(\text{resp. }\le)$.
By equation~(\ref{sig1}) we have the desired 
inequality $\nu_{\gamma}(F)<\nu_{\gamma}(E)$, 
$(\text{resp. }\le)$. 

\noindent$\Leftarrow$) 
We take an arbitrary 
exact sequence in $\A_{(\beta,t\omega)}$ 
\begin{equation}\label{ex}
0\to K\to E\to Q\to 0
\end{equation}
with 
$K\neq 0$ and $Q\neq 0$. 
We will show the inequality  
\begin{equation}\label{q}
\mu_{\sigma_{(\beta,t\omega)}}(\H^{-i}(Q)[i])
>\mu_{\sigma_{(\beta,t\omega)}}(E),
\ (\text{resp. }\ge)
\end{equation} 
if $\H^{-i}(Q)\neq0$ for $i=0,1$.
Then since $Z_{(\beta,t\omega)}(Q)
=Z_{(\beta,t\omega)}(\H^0(Q))
+Z_{(\beta,t\omega)}(\H^{-1}(Q)[1])$, 
we have the desired inequality 
$$\mu_{\sigma_{(\beta,t\omega)}}(Q)>
\mu_{\sigma_{(\beta,t\omega)}}(E),
\ (\text{resp. }\ge),$$ showing that
$E$ is $\sigma_{(\beta,
t\omega)}$-(semi)stable.

First we assume $\H^{-1}(Q)\neq0$ and 
show (\ref{q}). In fact we see that 
the inequality is always strict.
The fact that $E$ is 
a torsion free sheaf implies that $K$ is 
also a torsion free sheaf. Hence we have 
$\Im Z_{(\beta,t\omega)}(K)>0$. 
Since 
$$\Im Z_{(\beta,t\omega)}(E)=\Im Z_{(\beta,t\omega)}(K)
+\Im Z_{(\beta,t\omega)}(\H^0(Q))+\Im Z_{(\beta,t\omega)}
(\H^{-1}(Q)[1]),$$  
we see that $0\le\Im Z_{(\beta,t\omega)}(\H^{-1}(Q)[1])<
\Im Z_{(\beta,t\omega)}(E)=t\e$.  
The same argument as in the proof of 
Lemma~\ref{qlim}~(1) shows the strict inequality 
$\Re Z_{(\beta,t\omega)}\left(\H^{-1}(Q)[1]\right)<0$. 
Hence by the assumption that $\Re
Z_{(\beta,t\omega)}\left(E\right)\ge0$
we have the strict inequality 
$$\mu_{\sigma_{(\beta,t\omega)}}(E)
<\mu_{\sigma_{(\beta,t\omega)}}(\H^{-1}(Q)[1]).$$ 

Next we assume $\H^0(Q)\neq0$.
We take the cohomology 
long exact sequence of (\ref{ex})
in $\Coh(X)$; 
\begin{equation*}
0\to \H^{-1}(Q)\to K\to E\to \H^0(Q)\to 0.
\end{equation*}
We take $I:=\im(K\to E)$. 
Since the fact that $K,Q\in\A_{
(\beta,t\omega)}$ implies $\mu_{
t\omega}(K)>\mu_{t\omega}(\H^{-1}(Q))$,
we have $K\ncong \H^{-1}(Q)$. Hence
$I$ is not equal to $0$ and is torsion free.

If the strict inequality 
\begin{equation}\label{st}
\mu_{\omega}(I)<\mu_{\omega}(E)
\end{equation} 
holds we show 
a contradiction in the following way. 
We can write 
$$\mu_{t\omega}(E)-\mu_{t\omega}(I)=\frac
{t(r(I)c_1\cdot\omega-rc_1(I)\cdot\omega)}{rr(I)}.$$
By (\ref{st}) we have 
$(r(I)c_1\cdot\omega-rc_1(I)\cdot\omega)\in\Z_{>0}$.
Hence we get 
\begin{equation}
\label{ineq1}
\mu_{t\omega}(E)-\mu_{t\omega}(I)
\ge\frac{t}{r^2}.
\end{equation} 
On the other hand
since $K\to I$ is surjective, we
have the following inequalities 
\begin{equation*}
\beta\cdot t\omega< \mu_{t\omega\text{-min}}(K)\le 
\mu_{t\omega}(I).
\end{equation*}
Hence we get
\begin{equation}
\label{ineq2}
\mu_{t\omega}(E)-\mu_{t\omega}(I)
<\frac{c_1\cdot t\omega}{r}-\beta\cdot t\omega
=\frac{t\e}{r}
\end{equation}
by (\ref{ep}).
Combining (\ref{ineq1}) and (\ref{ineq2}) with the 
assumption that $\e\le\frac{1}{r}$, 
we get a contradiction.

In the case $r(I)=r$ and 
$\dim \H^0(Q)=1$ we have 
$\mu_{\omega}(I)<\mu_{\omega}(E)$. 
Hence we may assume that $0<\rk(I)<\rk(E)$ 
holds or that $\rk(I)=\rk(E)$ and 
$\dim(\H^0(Q))=0$ holds. 
In the latter case, we see that the slope 
$\mu_{\sigma_{(\beta,t\omega)}}(\H^0(Q))$  
is infinity and the desired inequality 
$\mu_{\sigma_{(\beta,t\omega)}}(E)< 
\mu_{\sigma_{(\beta,t\omega)}}(\H^0(Q))$ holds. 

We assume that $\rk(I)<\rk(E)$. 
Since $E$ is $(\gamma,
\omega)$-(semi)stable, 
\begin{equation*}
\left(\mu_{\omega}(E),\nu_{\gamma}(E)\right)<
\left(\mu_{\omega}(\H^0(Q)),\nu_{\gamma}
(\H^0(Q))\right),
\hspace{1cm}(\text{resp. }\le).
\end{equation*}
Then since $\mu_{\omega}(I)=\mu_{\omega}(E)$ by the
above argument, we have $$\mu_{\omega}(E)=
\mu_{\omega}(\H^0(Q))\hspace{5mm} \text{and} 
\hspace{5mm}\nu_{\gamma}(E)<\nu_{\gamma}(\H^0(Q)),
\ (\text{resp. }\le).$$ Hence by (\ref{sig1})
we get the desired 
inequality $\mu_{\sigma_{(\beta,t\omega)}}(E)
<\mu_{\sigma_{(\beta,t\omega)}}(\H^0(Q)),
\ (\text{resp. }\le)$.
\endproof

Here we assume that $\beta$ belongs to $\ns(X)\otimes
\mathbb{Q}$, or that $\gamma=\beta-\frac{1}{2}K_X$ is 
proportional to $\omega$ in $\ns(X)\otimes\R$.
In the latter case we have 
$\M_{X}(\ch(\alpha),\gamma,\omega)
=\M_{X}(\ch(\alpha),0,\omega)$
by (\ref{hil1}) and (\ref{hil2}). We recall that $\omega$
is an integral divisor.
Hence in both cases we have moduli spaces 
$M_{X}(\ch(\alpha),\gamma,\omega)$ of 
$\M_{X}(\ch(\alpha),\gamma,\omega)$
by \cite[Theorem~5.7]{MW}.

\begin{cor}\label{cor}Under the assumptions 
in the above theorem the moduli space 
$M_{X}(\ch(\alpha),\gamma,\omega)$ of 
$(\gamma,\omega)$-semistable sheaves  
corepresents the moduli functor 
$\M_{\D(X)}(\ch(\alpha),\sigma_{(\beta,t\omega)})$. 
In the case where $\gamma$ is 
proportional to $\omega$, or
$\gamma=0$, the open subset 
$M^s_{X}(\ch(\alpha),\omega)\subset M_{X}(\ch(\alpha),\omega)$ corepresents 
the functor $\M_{\D(X)}^s(\ch(\alpha),\sigma_{(\beta,t\omega)})$. 
\end{cor}
\proof
This follows directly from Theorem~\ref{limit}
and \cite[Lemma~3.31]{H}.
\endproof

By this corollary we get Theorem~\ref{geomintro}
in the introduction.

\section{Algebraic Bridgeland 
stability conditions}
\subsection{Moduli functors 
of representations of algebras}

For a finite dimensional 
$\C$-algebra $B$, we consider 
the abelian category mod-$B$ of 
finitely generated right 
$B$-modules and introduce the 
notion of $\theta_B$-stability 
of $B$-modules and families of 
$B$-modules over schemes 
following \cite{Ki1} . 

\begin{defn}\label{theta}
Let $\theta_B\colon K(B)\to \R$ 
be an additive function 
on the Grothendieck group $K(B)$. 
An object $N\in$ mod-$B$ is 
called $\theta_B$-semistable if 
$\theta_B(N)=0$ and every subobject 
$N'\subset N$ satisfies $\theta_B
(N')\ge 0$. Such an $N$ is called 
$\theta_B$-stable if the only 
subobjects $N'$ with $\theta_B(N')=0$ 
are $N$ and $0$.
\end{defn}

For $S\in(\text{scheme}/\C)$,  
define $\Coh_B(S)$ to be the 
category with objects $(F,\rho)$ 
for $F$ a coherent sheaf on $S$ 
and $\rho\colon B\to\Hom_{S}
(F,F)$ a $\C$-linear homomorphism 
with $\rho(ab)=\rho(b)\circ\rho(a)$ 
for each $a,b\in B$, and morphisms 
$\eta\colon(F,\rho)\to (F',\rho')$ 
to be morphisms of sheaves 
$\eta\colon F\to F'$ with $\eta\circ
\rho(a)=\rho'(a)\circ\eta$ in 
$\Hom_{S}(F,F')$ for all $a\in B$. 
It is easy to show $\Coh_B(S)$ is an 
abelian category. Let $\Vec_B(S)$ be 
the full subcategory of $\Coh_B(S)$ 
consisting of objects 
$(E,\rho)\in\Coh_B(S)$ where 
$E$ is locally free. 

\begin{defn}
Objects of $\Vec_B(S)$ are 
called families of 
$B$-modules over $S$ 
\emph{(\cite[Definition~5.1]
{Ki1})}.
\end{defn}

For $\alpha_B\in K(B)$ and an 
additive function 
$\theta_B\colon K(B)\to\R$ 
as in Definition~\ref{theta}, 
let $\M_{B}(\alpha_B,\theta_B)$ be 
the moduli functor which sends 
$S\in (\text{scheme}/\C)$ to 
the set $\M_{B}(\alpha_B,\theta_B)(S)$ 
consisting of isomorphism classes  
of families of $\theta_B$-semistable 
right $B$-modules $N$ with $[N]=
\alpha_B$.
Let $\M_B^s(\alpha_B,\theta_B)$ 
be the subfunctor of 
$\M_{B}(\alpha_B,\theta_B)$ 
corresponding to $\theta_B$-stable 
right $B$-modules.
There exist moduli 
spaces $M_B^s(\alpha_B,\theta_B)\subset
M_B(\alpha_B,\theta_B)$ of 
$\M_B^s(\alpha_B,\theta_B)$ 
and $\M_{B}(\alpha_B,\theta_B)$ 
\cite[Proposition~5.2]{Ki1}.

Here we recall the definition of 
the S-equivalence. Since any object
of mod-$B$ is finite dimensional
$\C$-vector space, any 
$\theta_B$-semistable $B$-module $N$
has a filtration, called \emph{
Jordan-H\"{o}lder filtration},
$$0=N_0\subset N_1\subset
\cdots \subset N_n=N$$ 
such that $N_i/N_{i-1}$ is 
$\theta_B$-stable for any $i$.
The \emph{grading} $Gr_{\theta_B}(N)
:=\oplus_iN_i/N_{i-1}$
does not depend on
a choice of a Jordan-H\"{o}lder 
filtration
up to isomorphism 
(for example, see \cite
[Proposition~1.5.2]{HL}). 
$\theta_B$-semistable $B$-modules 
$N$ and $N'$ are said to
be \emph{S-equivalent} if $Gr_{\theta_B}(N)
\cong Gr_{\theta_B}(N')$.

\begin{prop}
\emph{\bf{(cf. \cite[Proposition~3.2]{Ki1})}}\label{S-eq}
For $B$-modules $N$ and $N'$ with 
$[N]=[N']=\alpha_B\in K(B)$, 
$N$ and $N'$ define the same point
of $M_B(\alpha_B,\theta_B)$ if and only if
they are S-equivalent to each other.
\end{prop}

\subsection{Algebraic Bridgeland 
stability conditions}
\label{algbri}
Let $X$ be a smooth projective 
surface.
An object $E\in\D(X)$ is said to be
\emph{exceptional} if
$$\Hom^k_{\D(X)}(E,E)=
\begin{cases}
\C&\text{if }k=0\\
0&\text{otherwise}.
\end{cases}$$
An \emph{exceptional
collection} in $\D(X)$ 
is a sequence of
exceptional objects 
$\mathfrak{E}
=\left(E_0,\cdots , E_n\right)$
of $\D(X)$ such that 
$$n\ge i>j\ge 1\implies
\Hom^k_{\D(X)}(E_i,E_j)=0
\text{ for all }k\in\Z.$$
The exceptional collection
$\mathfrak{E}$ is said to be
full if 
$E_0,\cdots,E_n$ generates 
$\D(X)$, namely the smallest
full triangulated subcategory
containing $E_0,\cdots,E_n$
coincides with $\D(X)$. 
The exceptional collection
$\mathfrak{E}$ is said to
be strong 
if for all $1\le i,j\le n$
one has 
$$\Hom^k_{\D(X)}(E_i,E_j)=0
\text{ for } k\neq0.$$

We assume that $\D(X)$ has 
a full strong 
exceptional collection 
$\mathfrak{E}=\left(E_0,
\cdots,E_n\right)$ on $\D(X)$. 
We put $\E:=E_0\oplus\cdots
\oplus E_n$, $B_{\E}:=\End_X(\E)$. 
By Bondal's theorem \cite{Bo}
we have an equivalence 
$$\Phi_{\E}\colon \D(X)\cong\D(B_{\E})
\colon E\mapsto\mathbf R\Hom_X(\E,E).$$ 
We obtain the heart $\A_{\E}
\subset\D(X)$ by pulling back 
mod-$B_{\E}$ via the equivalence 
$\Phi_{\E}$. 
The equivalence $\Phi_{\E}$ 
induces an isomorphism $\varphi_{\E}
\colon K(X)\cong K( B_\E)$ of the 
Grothendieck groups. 

For a stability function $Z$ on $\A_\E$
and $\alpha\in K(X)$, we define  
$\theta^{\alpha}_Z\colon K(B_\E)\to\R$ by  
\begin{equation}\label{det}
\theta^{\alpha}_{Z}(\beta):=
\begin{vmatrix}
\Re Z(\varphi_\E^{-1}(\beta)) & \Re Z(\alpha)\\
\Im Z(\varphi_\E^{-1}(\beta)) & \Im Z(\alpha)
\end{vmatrix}
\end{equation}
for any $\beta\in K(B_\E)$. 
Then for an object 
$E\in\A_\E$ with 
$[E]=\alpha \in K(X)$, 
$E$ is $Z$-(semi)stable
if and only if $\Phi_\E(E)$ is 
$\theta^{\alpha}_Z$-(semi)stable. We also notice that 
by the existence of full exceptional collection, 
$K(X)$ is isomorphic to
the numerical Grothendieck group $K(X)/
K(X)^{\perp}$. Hence for $E\in\D(X)$ 
the class $[E]$ is equal to $\alpha$ 
in $K(X)$ if and only if $\ch(E)=\ch(\alpha)$.

\begin{prop}\label{algebraic} 
The moduli space $M_{B_\E}(\varphi_\E(\alpha),
\theta^{\alpha}_Z)$ $($resp. $M^s_{B_\E}(\varphi_\E(\alpha),
\theta^{\alpha}_Z)$ $)$ corepresents the moduli functor 
$\M_{\D(X)}(\ch(\alpha),\sigma)$ 
$($resp. $\M^s_{\D(X)}(\ch(\alpha),
\sigma)$ $)$ for any $\alpha\in 
K(X)$, $\sigma=(Z,\A_{\E})\in\Stab(\A_{\E})$.
\end{prop}
\proof
We only give the proof for the moduli functor
$\M_{\D(X)}(\ch(\alpha),\sigma)$,
since a similar argument also holds 
for the other moduli functor
$\M^s_{\D(X)}(\ch(\alpha),\sigma)$
corresponding to stable objects.
We show that 
\begin{equation}
\label{shiso}
{}^{sh}\M_{\D(X)}
(\ch(\alpha),\sigma)\cong{}^{sh}
\M_{B_\E}(\varphi_{\E}(\alpha),\theta^{\alpha}_Z). 
\end{equation}
Then, since $M_{B_\E}(\varphi_{\E}(\alpha),\theta^{\alpha}_Z)$
corepresents $\M_{B_\E}(\varphi_{\E}(\alpha),
\theta^{\alpha}_Z)$, the assertion holds
by (\ref{sh}). By the remark after
(\ref{shi}), 
to establish (\ref{shiso})
it is enough to give a functorial
isomorphism
\begin{equation}\label{map}
\M_{\D(X)}
(\ch(\alpha),\sigma)(S)\cong\M_{B_\E}(\varphi_{\E}(\alpha),
\theta^{\alpha}_Z)(S),
\end{equation}
for every affine scheme
$S=\Spec R$. 
We consider $X_S:=X\times S$, 
projections $p$ and $q$ from
$X_S$ to $X$ and $S$, the pull 
back $\E_S:=p^{\ast}\E$ 
of $\E$ and $R$-algebra
$B_{\E_S}:=\Hom_{X_S}
(\E_S,\E_S)$. 
Since $B_{\E_S}\cong R\otimes {B_\E}$,
we have mod-$B_{\E_S}\cong\Coh_{B_\E}(S)$.
From \cite[Lemma~8]{TU} we see that 
via the above identification 
$\Phi_{\E_S}(\ \cdot\ ):=\mathbf R
\Hom_{X_S}(\E_S,\ \cdot\ )$
gives equivalences  
$$\D(X_S)\cong\D(\Coh_{B_\E}(S)),
\hspace{1.5cm} 
D^-(X_S)\cong D^-(\Coh_{B_\E}(S)).$$
These equivalences are 
compatible with pull backs, 
that is, the following diagram 
is commutative
$$\xymatrix{D^-(X_S)
\ar[d]_{\mathbf L f^{\ast}}
\ar[rr]^{\Phi_{\E_S}}&&
D^-(\Coh_{B_\E}(S))
\ar[d]^{\mathbf L f^{\ast}}
\ar@{}[lld]|\circlearrowleft\\
D^-(X_{S'})
\ar[rr]^{\Phi_{\E_{S'}}}
&&D^-(\Coh_{B_\E}(S'))}
$$ 
for every morphism $f\colon 
S'\to S$ of affine schemes.
In the following we show that 
this equivalence $\Phi_{\E_S}$ 
defines an 
isomorphism (\ref{map}).
  
For any $S$-valued point 
$E$ of $\M_{\D(X)}(\ch(\alpha),\sigma)$, 
by the above diagram the fact that
$E\in\M_{\D(X)}(\ch(\alpha),\sigma)(S)$ 
implies that
$\mathbf L \iota_s^{\ast}\Phi_{\E_S}(E)
\in D^-(\Coh_{{B_\E}}(\{s\}))\cong
D^-({B_\E})$ is a 
$\theta_Z^{\alpha}$-semistable 
${B_\E}$-module for any closed point $s\in S$, 
where $\iota_s\colon \{s\}\to S$ is 
the embedding.  
By the standard argument using 
the spectral sequence 
(for example, \cite[Lemma 3.31]{H}), 
we see that $\Phi_{\E_S}(E)$ belongs to 
$\Vec_{B_\E}(S)\subset\Coh_{B_\E}(S)$.
Hence $\Phi_{\E_S}$ defines a
desired map. 
We see that this map
is an isomorphism since 
$\Phi_{\E_S}^{-1}$ gives
the inverse map by a similar 
argument.
\endproof
By this proposition we get Proposition~\ref{algintro}
in the introduction.
\begin{defn}\label{alg}
$\sigma\in\Stab(X)$ is called an 
algebraic Bridgeland stability 
condition associated to 
the full strong exceptional 
collection $\mathfrak{E}=
\left(E_0,\ldots,E_n\right)$
if $\sigma$ is 
$\grp$-equivalent to $(Z,\A_\E)$ 
for some $Z\colon K(X)\to\C$, where
$\E=E_0\oplus\cdots\oplus E_n$.
\end{defn}

\subsection{Full strong exceptional
collections on $\PP^2$}
\label{fsecp2}
In the rest of the paper, we assume 
that $X = \PP^2$ and $H$ is the 
hyperplane class on $\PP^2$. 
We put $\mo_{\PP^2}(1)
:=\mo_{\PP^2}(H)$ and denote 
the homogeneous coordinates of
$\PP^2$ by $[z_0:z_1:z_2]$.
We introduce two types of full
strong exceptional collections
$\mathfrak{E}_k$ and $\mathfrak{
E}'_k$ on $\PP^2$ for each $k\in\Z$
as follows, 
$$\mathfrak{E}_k:=\left(\mo_{\PP^2}(k+1),
\Omega^1_{\PP^2}(k+3),\mo_{\PP^2}(k+2)
\right),\hspace{15mm}
\mathfrak{E}'_k:=\left(\mo_{\PP^2}(k),
\mo_{\PP^2}(k+1),\mo_{\PP^2}(k+2)\right).
$$
We put 
$$\mathcal{E}_k:=\mo_{\PP^2}(k+1)
\oplus\Omega^1_{\PP^2}(k+3)
\oplus\mo_{\PP^2}(k+2),\hspace{1.5cm} 
\E'_k:=\mo_{\PP^2}(k)\oplus\mo_{\PP^2}(k+1)
\oplus\mo_{\PP^2}(k+2)$$ and
$B:=\End_{\PP^2}(\mathcal{E}_k)$, 
$B':=\End_{\PP^2}(\E'_k)$,
which do not depend on $k$ up to
natural isomorphism.
Using the notation in 
\S~\ref{algbri}, we define 
functors 
$$\Phi_k:=\Phi_{\E_k}\colon 
\D(\PP^2) \cong \D(B), \ \ \ 
\Phi'_k:=\Phi_{\mathcal{E}'_k}\colon 
\D(\PP^2) \cong \D(B'),$$
induced isomorphisms
$\varphi_k:=\varphi_{\mathcal{E}_k}
\colon K(\PP^2)\cong K(B)$, 
$\varphi'_k:=\varphi_{\E'_k}
\colon K(\PP^2)\cong K(B')$
and full subcategories $\A_k:=
\A_{\mathcal{E}_k}$, $\A'_k
:=\A_{\E'_k}$ of $\D(\PP^2)$.

To explain finite dimensional 
algebras $B$ and $B'$ 
we introduce some notations. 
For any $l\in\Z$, we denote by $z_i$ the morphism 
$\mo_{\PP^2}(l)\to\mo_{\PP^2}(l+1)$
defined by multiplication of 
$z_i$ for $i=0,1,2$. We put
$V:=\C e_0\oplus\C e_1\oplus\C e_2$ and 
denote $i$-th projection and $i$-th embedding 
by $e_i^{\ast}\colon V\to\C$ and 
$e_i\colon \C\to V$ for $i=0,1,2$.
We consider the exact sequence for each $k\in\Z$
\begin{equation}\label{euler}
0\to\Omega^1_{\PP^2}(k+3)
\xrightarrow{\iota}\mo_{\PP^2}(k+2)\otimes V
\xrightarrow{j}\mo_{\PP^2}(k+3)\to 0,
\end{equation}
where we put $j:=z_0\otimes e_0^{\ast}+
z_1\otimes e_1^{\ast}+
z_2\otimes e_2^{\ast}$
and identify $\Omega^1_
{\PP^2}(k+3)$ with $\ker j$.
We define morphisms 
$p_i\colon\Omega^1_{\PP^2}(k+3)
\to\mo_{\PP^2}(k+2)$ by 
$p_i:=(\id_{\mo_{\PP^2}(k+2)}
\otimes e_i^{\ast})\circ\iota$ and
$q_i\colon\mo_{\PP^2}(k+1)
\to\Omega^1_{\PP^2}(k+3)$
by $q_i:=z_{i+2}\otimes e_{i+1}-
z_{i+1}\otimes e_{i+2}$ for 
$i\in\Z/3\Z$.

We introduce the following
quiver $Q$ with
$3$ vertices $\{v_0,v_1,
v_2\}$ and $6$ arrows
$\{\gamma_0,\gamma_1,\gamma_2,
\delta_0,\delta_1,\delta_2\}$
\begin{equation*}
\stackrel{v_0}{\bullet} 
\xleftarrow{\hspace{0.8cm}
\gamma_i\hspace{0.8cm}} 
\stackrel{v_1}{\bullet} 
\xleftarrow{\hspace{0.8cm}
\delta_j\hspace{0.8cm}} 
\stackrel{v_2}{\bullet} 
\hspace{2cm} (i,j=0,1,2)
\end{equation*}
and consider ideals $J$ 
and $J'$ of the path 
algebra $\C Q$ defined 
as follows. $J$ and $J'$ 
are two-sided ideals
generated by $\{\gamma_i\delta_j+
\gamma_j\delta_i\mid i,j=0,1,2\}$
and $\{\gamma_i\delta_j-
\gamma_j\delta_i
\mid i,j=0,1,2\}$, respectively.
We have isomorphisms 
\begin{equation}\label{final}
\rho\colon\C Q/J\cong B\colon
\gamma_i,\delta_j\mapsto p_i, q_j,
\hspace{0.7cm}\rho'\colon\C Q/J'\cong B'\colon
\gamma_i,\delta_j\mapsto z_i,z_j.
\end{equation} 
These isomorphisms $\rho$ and $\rho'$ map 
vertices $v_0, v_1, v_2\in \C Q/J$ (resp. $\C Q/J'$)
to idempotent elements $$\rho(v_0)=\id_{\mo_{\PP^2}(k+2)},\ \ 
\rho(v_1)=\id_{\Omega^1_{\PP^2}(k+3)},\ \ \rho(v_2)=\id_{\mo_{\PP^2}(k+1)}\in B$$ 
$$\left(\text{resp. }\rho'(v_0)=\id_{\mo_{\PP^2}(k+2)},\ \ \rho'(v_1)=\id_{\mo_{\PP^2}(k+1)},\ \ 
\rho'(v_2)=\id_{\mo_{\PP^2}(k)}\in B'\right).$$ 
They also map
$\gamma_i, \delta_j\in \C Q/J$ 
(resp. $\C Q/J'$) to 
$$\rho(\gamma_i)=p_i,\ \ 
\rho(\delta_j)=q_j\in B\ \ \ \ 
\left(\text{resp. }\rho'(\gamma_i)=z_i,\ \ 
\rho'(\delta_j)=z_j\in B'\right)$$ for $i,j=0,1,2$.
We identify $B$ and $B'$ with 
$\C Q/J$ and $\C Q/J'$ via isomorphisms
$\rho$ and $\rho'$. 

For any finitely generated 
right $B$-module $N$, we
consider the right action 
on $N$ of a path $p$ of $Q$ 
as a pull back by $p$ and 
denote it by $p^*$. Notice 
that vertices $v_i's$ are 
regarded as paths 
with the length $0$. 
We have 
the decomposition $N=Nv^{\ast}_0\oplus 
Nv^{\ast}_1\oplus Nv^{\ast}_2$ as a vector 
space. This gives the dimension 
vector $\dimv(N)=(\dim_{\C} Nv^{\ast}_0,
\dim_{\C} Nv^{\ast}_1,\dim_{\C} Nv^{\ast}_2)$ 
of $N$ and an isomorphism $\dimv
\colon K(B)\cong\Z^{\oplus 3}$. 
The $B$-module structure of $N$ is 
written as;
$$Nv^{\ast}_0\xrightarrow{\gamma^{\ast}_i}
Nv^{\ast}_1\xrightarrow{\delta^{\ast}_j}
Nv^{\ast}_2 \ \ (i,j=0,1,2).$$
We sometimes use notation 
$\gamma^{\ast}_i|_N$ and $\delta^
{\ast}_j|_N$ to avoid confusion.
We define 
$B$-modules $\C v_i$ for 
$i=0,1,2$ as follows. As vector
spaces $\C v_i=\C$ and can be decomposed
by $(\C v_i)v^{\ast}_i=\C$, $(\C v_i)v^{\ast}_j=0$ 
for $j\neq i$. Actions of $B$ are defined
in obvious way. They are simple objects
of mod-$B$ and we have
\begin{equation}
\text{mod-}B=\langle\C v_0, \C v_1,
\C v_2\rangle
\end{equation}
as a full subcategory
of $\D(B)$.
Similar results hold 
for $B'$ and we use similar notations 
for $B'$. 

Since $\mo_{\PP^2}(k-1)[2]$,
$\mo_{\PP2}(k)[1]$ and $\mo_{\PP^2}(k+1)$
correspond to $B$-modules
$\C v_0$, $\C v_1$ and 
$\C v_2$ via $\Phi_k$, we have 
\begin{equation*}
\A_k=\langle\mo_{\PP^2}(k-1)[2],
\mo_{\PP2}(k)[1],\mo_{\PP^2}(k+1)
\rangle.
\end{equation*}
Similarly we have 
$$\A'_k=\langle \mo_{\PP^2}(k-1)[2],
\Omega^1_{\PP2}(k+1)[1],
\mo_{\PP^2}(k)\rangle.
$$
On the other hand,
$\mo_{\PP^2}(k+1)$, $\Omega^1_{
\PP^2}(k+3)$ and $\mo_{\PP^2}(k+2)$
correspond to $B$-modules
$B$, $v_1B$ and $v_2B$ 
via $\Phi_{k}$.
Similarly $\mo_{\PP^2}(k)$,
$\mo_{\PP^2}(k+1)$ and 
$\mo_{\PP^2}(k+2)$ correspond to 
$B'$-modules $B'$, 
$v_1B'$ and $v_2B'$
via $\Phi'_{k}$.
They are projective modules
and we can compute Ext groups
by using them. Hence we get 
the following lemma.

\begin{lem}\label{bei1}
For bounded complexes $E, F$ 
of coherent 
sheaves on $\PP^2$, the following 
hold for each $k\in\Z$.\\
$(1)$ By $E^i$, we denote
each term of complex $E$. 
We assume that $($\emph{i}$)$ 
$E^i$ is a direct sum 
of $\mo_{\PP^2}(k+1)$, $\Omega^1_
{\PP^2}(k+3)$ and $\mo_{\PP^2}(k+2)$ 
for any $i\in\Z$ and $F$ belongs to 
$\A_k$, or that $($\emph{ii}$)$ 
$E^i$ is a direct sum 
of $\mo_{\PP^2}(k)$, $\mo_
{\PP^2}(k+1)$ and $\mo_{\PP^2}(k+2)$
for any $i\in\Z$ and $F$ belongs 
to $\A'_k$. Then
the complex $\mathbf R \Hom_{\PP^2}(E,F)$ 
is quasi-isomorphic to the 
following complex
\begin{equation}\label{rhom}
\cdots \to\Hom_{\D(\PP^2)}
(E^{-i},F)\xrightarrow{d^i}
\Hom_{\D(\PP^2)}(E^{-i-1},F)
\to\cdots,
\end{equation}
where $\Hom_{\D(\PP^2)}
(E^{-i},F)$ lies on degree $i$
and $d^i$ is defined by $$d^i
(f):=f\circ d^{-i-1}_{E}\colon
E^{-i-1}\to F\hspace{.5cm}
\text{ for }f\in\Hom_{\PP^2}(E^{-i},F).$$
In particular, we have 
$\Hom_{\D(\PP^2)}(E,F[i])\cong
\ker d^i/\im d^{i-1}.$
\\
$(2)$ If $E$ belongs to $\A_k$
$($resp. $\A'_k)$, then
we have the following isomorphism 
in $\D(\PP^2)$
$$E\cong\left(\mo_{\PP^2}(k-1)^{\oplus
a_0}\to\mo_{\PP^2}(k)^{\oplus
a_1}\to\mo_{\PP^2}(k+1)^{\oplus
a_2}\right),$$
$$\left(\text{resp. }E\cong
\left(\mo_{\PP^2}(k-1)^{\oplus
a_0}\to\Omega^1_{\PP^2}(k+1)^{\oplus
a_1}\to\mo_{\PP^2}(k)^{\oplus
a_2}\right)\right),$$
where $(a_0,a_1,a_2)\in\Z^3_{\ge0}$ and 
$\mo_{\PP^2}(k+1)^{\oplus
a_2}$ $($resp. $\mo_{\PP^2}(k)
^{\oplus a_2})$ lies on degree $0$.\\
\end{lem}
\proof
$(1)$ We only prove (i). We put 
$N:=\Phi_k(E)$, $M:=\Phi_k(F)$. 
Then by the 
assumption the each term $N^i$ of the complex $N$ is a 
direct sum of $B$, $v_1B$ and $v_2B$ for any $i$.
Hence $N^i$ is a projective module. Furthermore
since the fact $F\in\A_k$ implies that
$M$ is a $B$-module, 
$\mathbf R\Hom_{\PP^2}(E,F)\cong\mathbf R
\Hom_{B}(N,M)$ is quasi-isomorphic to 
the following complex
\begin{equation*}
\cdots \to\Hom_{B}
(N^{-i},M)\xrightarrow{d^i}
\Hom_{B}(N^{-i-1},M)
\to\cdots.
\end{equation*}
Via $\Phi_k$ this complex coincides with
(\ref{rhom}).\\
$(2)$ For any object $E\in\A_{k}$
we consider the $B$-module $N=\Phi_k(E)$.
If we put $\dimv(N)=(a_0,a_1,a_2)$,
then $N$ can be obtained by extensions 
\begin{equation}\label{ext1}
0\to(\C v_1)^{\oplus a_1}\to N'\to
(\C v_0)^{\oplus a_0}\to0,
\end{equation}
\begin{equation}\label{ext2}
0\to (\C v_2)^{\oplus a_2}\to N\to N'
\to 0.
\end{equation}
Since $\Phi_k(\mo_{\PP^2}(k-1)[1])=\C v_0[-1]$
and $\Phi_k(\mo_{\PP^2}(k)[1])=\C v_1$, we have a 
homomorphism $$f\colon \mo_{\PP^2}(k-1)^{\oplus a_0}\to
\mo_{\PP^2}(k)^{\oplus a_1}$$ in $\Coh(\PP^2)$
such that $\Phi_k(C(f)[1])\cong N'$,
where $C(f)$ is the mapping cone of $f$.
From (\ref{ext2}) $E$ can be 
obtained as a mapping cone 
of a certain homomorphism in $\Hom
_{\D(\PP^2)}(C(f),\mo_{\PP^2}(k+1)
^{\oplus a_2})$, since $\Phi_k(\mo_{
\PP^2}(k+1))=\C v_2$. By (1) 
this homomorphism is identified
with a homomorphism $$g\colon
\mo_{\PP^2}(k)^{\oplus a_1}\to\mo_{\PP^2}
(k+1)^{\oplus a_2}$$ in $\Coh(\PP^2)$ satisfying $g\circ f=0$.
Thus $E$ is isomorphic to the following
complex 
$$\left(\mo_{\PP^2}(k-1)^{\oplus a_0}\stackrel{f}
{\to}\mo_{\PP^2}(k)^{\oplus a_1}\stackrel{g}
{\to}\mo_{\PP^2}(k+1)^{\oplus a_2}\right),$$
where $\mo_{\PP^2}(k+1)^{\oplus a_2}$ lies
on degree $0$.
\endproof
The vector $(a_0,a_1,a_2)\in\Z_{\ge 0}^3$
in Lemma~\ref{bei1}~(2) coincides with 
$\dimv(\Phi_k(E))$ and is explicitly computed
from $\ch(E)=(r,sH,\ch_2)$.
For example, we assume that $E$ belongs 
to $\A_1$. Since 
\begin{equation}\label{dim}
\ch(\mo_{\PP^2}[2])=(1,0,0),\hspace{5mm}
\ch(\mo_{\PP^2}(1)[1])=-(1,H,\frac{1}{2}),
\hspace{5mm}\ch(\mo_{\PP^2}(2))=(1,2H,2),
\end{equation}
we have $(a_0,a_1,a_2)
=r(1,0,0)-\frac{s}{2}(3,4,1)+\ch_2(1,2,1)$.

\section{Proof of Main Theorem~\ref{inmain}}
\label{quiver}

In this section we fix $\alpha
\in K(\PP^2)$ with $\ch(\alpha)=
(r,sH,\ch_2)$ and $0<s\le r$.
In the sequel, we sometimes identify $\ns(\PP^2)$ with 
$\Z$ by the isomorphism
$\ns(\PP^2)\cong \Z\colon
\beta\mapsto\beta\cdot H$.

\subsection{Wall-and-chamber structure}
\label{main}
We consider the full strong 
exceptional collection 
$\mathfrak{E}_1=\left(\mo_{\PP^2}(2), 
\Omega_{\PP^2}^1(4), 
\mo_{\PP^2}(3)\right)$ on $\PP^2$,
the equivalence $\Phi_1(\ \cdot\ )=\mathbf R
\Hom_{\PP^2}(\E_1,\ \cdot\ )\colon
\D(\PP^2)\cong\D(B)$
and the induced isomorphism $\varphi_1
\colon K(\PP^2)\cong K(B)$,
where $\E_1=\mo_{\PP^2}(2)\oplus 
\Omega_{\PP^2}^1(4)\oplus 
\mo_{\PP^2}(3)$ and 
$B=\End_{\PP^2}(\E_1)$.
We consider the plane 
$\varphi_1(\alpha)^{\perp}:=
\{\theta_1\in\Hom_{\Z}(K({B}),\R)
\mid\theta_1(\varphi_1(\alpha))=0\}$ 
and define a subset $W_1\subset 
\varphi_1(\alpha)^{\perp}$ as follows. 
A subset $W_1$ consists of 
elements $\theta_1\in
\varphi_1(\alpha)^{\perp}$ satisfying
that there exists a 
$\theta_1$-semistable
${B}$-module $N$ with $[N]=
\varphi_1(\alpha)$ such that 
$N$ has a proper nonzero submodule 
$N'\subset N$ with 
$\theta_1(N')=0$ 
and $[N']\notin\Q_{>0}\varphi_1(\alpha)$ 
in $K(B)$. 
The subset $W_1$ is a union of 
finitely many rays in 
$\varphi_1(\alpha)^{\perp}$.
These rays are called walls and 
the connected components of 
$\varphi_1(\alpha)^{\perp}
\setminus W_1$ are called 
chambers. 

We take a line $l_1$ 
in $\varphi_1(\alpha)^{\perp}$
defined by $l_1:=\{\theta_1
\in\varphi_1(\alpha)^{\perp}
\mid\theta_1(\varphi_1(\mo_x))=0\}$, where 
$\mo_x$ is the structure sheaf 
of a point $x\in\PP^2$. 
We take a chamber 
$C^{\PP^2}_{\varphi_1(\alpha)}\subset\varphi_1(\alpha)^{\perp}$, 
if any, such that the closure 
intersects with 
$l_1$ and there exists an
element $\theta_1\in C^{\PP^2}_{\varphi_1(\alpha)}$
satisfying the inequality 
$\theta_1(\varphi_1(\mo_x))>0$ and 
$M_{B}(-\varphi_1(\alpha),
\theta_1)\neq\emptyset$.
These conditions characterize 
$C^{\PP^2}_{\varphi_1(\alpha)}$ uniquely.

We have the following theorem,
which gives a proof of (i) in
Main Theorem~\ref{inmain}. 
The proof of Theorem~\ref{our} in the 
next subsection shows that 
if there is not such a chamber
$C^{\PP^2}_{\varphi_1(\alpha)}\subset\varphi_1(\alpha)^{\perp}$, 
then 
$M_{\PP^2}(\ch(\alpha),H)=\emptyset$.

\begin{thm}\label{our}
The map $E\mapsto\Phi_1(E[1])$ gives 
an isomorphism 
$$M_{\PP^2}(\ch(\alpha),H)\cong 
M_{{B}}(-\varphi_1(\alpha),\theta_1)$$ 
for any $\theta_1\in C^{\PP^2}_{\varphi_1(\alpha)}$.
This isomorphism keeps open subsets 
consisting of stable objects.
\end{thm}

Here we remark that if we assume 
$M_{\PP^2}(\ch(\alpha),H)\neq\emptyset$,
then $\dim M_{\PP^2}(\ch(\alpha),H)=s^2
-r^2+1-2r\ch_2\ge 0$.
Hence we have $\ch_2\le\frac{1}{2}$. 
We see that $\ch_2=\frac{1}{2}$ if and 
only if $\ch(\alpha)=(1,1,\frac{1}
{2})$.

\subsection{Proof of Theorem~\ref{our}}
\label{proof}
We will 
find Bridgeland stability 
conditions $\sigma$ in 
$\Stab(\A_1)\cap \{\sigma_{(bH,
tH)}\in\Stab(\PP^2)\mid t >0\}\grp$ for 
suitable $b\in\R$  
and obtain Theorem~\ref{our}. 

We put $\mathbf H=\big\{r\exp
(\sqrt{-1}\pi\phi)\mid r>0 
\text{ and }0<\phi\le 1\big\}$ 
the strict upper half-plane and 
$F_0=\mo_{\PP^2}[2]$, $F_1=
\mo_{\PP^2}(1)[1]$
and $F_2=\mo_{\PP^2}(2)$. 
The full subcategory 
$\A_1$ of $\D(\PP^2)$ is generated by 
$F_0$, $F_1$ and $F_2$, 
\begin{equation}\label{heart}
\A_1=\langle \mo_{\PP^2}[2],
\mo_{\PP^2}
(1)[1], \mo_{\PP^2}(2)
\rangle.
\end{equation}

Since $K(\PP^2)=\Z [F_0]\oplus
\Z [F_1]\oplus\Z [F_2]$, 
a stability function $Z$ on 
$\A_1$ is identified with the 
element $(Z(F_0),Z(F_1),Z(F_2))$ 
of $\mathbf H^3$.
Furthermore since the category 
$\A_1\cong\Mod$-$B$ has finite length, 
all stability functions on 
$\A_1$ satisfy the 
Harder-Narasimhan property. 
Hence $\Stab(\A_1)\cong
\mathbf H^3$. 

For $\sigma=(Z,\A_1)\in
\Stab(\A_1)$, we put $Z(F_i)=x_i+
\sqrt{-1}y_i\in\mathbf H^3$
and consider the conditions 
for $\sigma$ to be geometric. 
In the next lemmas we consider  
the condition~1 of 
Proposition~\ref{geometric1}. 
For any point $x\in\PP^2$ we take a 
resolution of $\mo_x$
\begin{equation}\label{skys}
0\to\mo_{\PP^2}\to\mo_{\PP^2}(1)^{\oplus 2}
\to\mo_{\PP^2}(2)\to\mo_x\to 0.
\end{equation}
Hence from (\ref{heart}) we have 
$\mo_x\in\A_1$ and $[\mo_x]=[F_0]+
2[F_1]+[F_2]\in K(\PP^2)$.
\begin{lem}\label{geom1}
For any subobject $E$ of $\mo_x$ in $\A_1$, 
the class $[E]$ in $K(\PP^2)$ is equal to 
$[F_2]$, $[F_1]+[F_2]$ or $2[F_1]+[F_2]$.
\end{lem}
\proof
If the conclusion is not true,
we can find a subobject
$\mathcal{F}[i]\subset\mo_x$ in $\A_1$ with 
$\mathcal{F}$ a nonzero sheaf on $\PP^2$ 
and $i= 1$ or $2$; 
for example, if $E$ is a subobject of $\mo_x$
in $\A_1$ and $[E]=[F_0]+[F_1]+[F_2]$
in $K(\PP^2)$, then by 
Lemma~\ref{bei1}~(2), $E$ is written as 
$$E=\left(\mo_{\PP^2}\stackrel{f}{\to}\mo_{\PP^2}(1)
\stackrel{g}{\to}\mo_{\PP^2}(2)\right).$$ 
If $g=0$ and $f\neq0$, then $E=\mo_{\ell}(1)[1]
\oplus\mo_{\PP^2}(2)$, 
where $\ell$ is a line on $\PP^2$ determined by 
$\mo_{\ell}(1)=\coker f$. If $g=f=0$, then $E=\mo_{\PP^2}[2]
\oplus\mo_{\PP^2}(1)[1]\oplus\mo_{\PP^2}(2)$.
If $g\neq 0$, then we have a distinguished
triangle 
$$\mo_{\ell'}(2)\to E\to\mo_{\PP^2}[2]\to
\mo_{\ell'}(2)[1]$$
for a line $\ell'$ on $\PP^2$ determined by 
$\mo_{\ell'}(2)=\coker g$. The fact that
$\Hom_{\D(\PP^2)}(\mo_{\PP^2}[2],\mo_{\ell'}(2)[1])=0$
implies $E=\mo_{\PP^2}[2]\oplus\mo_{\ell'}(2)$.

However the fact that 
$\Hom_{\D(\PP^2)}(\mathcal{F}[i],\mo_x)=0$ 
for $i\ge 1$ contradicts the fact that 
$\mathcal{F}[i]$ is a nonzero 
subobject of $\mo_x$ in $\A_1$. 
\endproof

\begin{lem}\label{geom2}
For $\sigma=(Z,\A_1)\in
\Stab(\A_1)$, $\mo_x$ is 
$\sigma$-stable for each 
$x\in\PP^2$ if and only if 
$(a)$, $(b)$ and $(c)$ hold;
$$(a)\begin{vmatrix}
x_2 & x_0+2x_1+x_2\\
y_2 & y_0+2y_1+y_2
\end{vmatrix}
>0, \ \ 
(b)\begin{vmatrix}
x_1+x_2 & x_0+2x_1+x_2\\
y_1+y_2 & y_0+2y_1+y_2
\end{vmatrix}
>0, \ \ 
(c)\begin{vmatrix}
2x_1+x_2 & x_0+2x_1+x_2\\
2y_1+y_2 & y_0+2y_1+y_2
\end{vmatrix}
>0.
$$
\end{lem}
\proof
By lemma~\ref{geom1}, it is enough to show 
$\phi(\beta)<\phi(\mo_x)$ for each $\beta=[F_2]$, 
$[F_1]+[F_2]$, $2[F_1]+[F_2]$, 
where $\phi(\beta)$ is the 
phase of $Z(\beta)\in\C$.
It is equivalent to 
$$\begin{vmatrix}\Re Z(\beta) 
& \Re Z(\mo_x) \\
\Im Z(\beta) & \Im Z(\mo_x)
\end{vmatrix} > 0,$$
which is equivalent to (a), (b) and (c)
for the case $\beta=[F_2]$, $[F_1]+[F_2]$ and 
$2[F_1]+[F_2]$ respectively. Hence the assertion follows.   
\endproof

By Lemma~\ref{geom2} and some 
easy calculations, we can find 
Bridgeland stability conditions 
$\sigma^b=(Z^b,\A_1)$ with $0<b<1$
which satisfy 
the conditions~1 and 2 in 
Proposition~\ref{geometric1}
as follows. 
We put $x_0:=-b$, $x_1:=-1+b$, 
$x_2:=-3b+3$ and $y_0=y_1=0$, 
$y_2=1$, that is,  
\begin{equation}\label{def1}
Z^b(F_0):=-b, \ \ Z^b(F_1):=-1+b, \ \ 
Z^b(F_2):=-3b+3+\sqrt{-1}.
\end{equation}
${\sigma}^{b}=(Z^{b},\A_1)
\in\Stab(\PP^2)$  satisfies the conditions 
$(a)$, $(b)$ and $(c)$ in Lemma~\ref{geom2}. 
The vector $\pi({\sigma}^{b})$ is written as 
$$\pi({\sigma}^{b})=u+\sqrt{-1}v\in
\N(\PP^2)\otimes\C$$ 
with $u=(2b-1,(b+\frac{1}{2})H,b),v=(-1,
-\frac{1}{2}H,0)\in\N(\PP^2)$. If we put 
$$T^{-1}:=
\begin{pmatrix}
b-\frac{1}{2}&2b^2-2b-
\frac{1}{2}\\
\sqrt{b-b^2}&(2b-1)\sqrt{b-b^2}
\end{pmatrix}\in
\GL^+(2,\R),
$$
then $\pi({\sigma}^{b})T
=\exp(bH+\sqrt{-1}
\sqrt{b-b^2}H);$
$$\begin{pmatrix}b-\frac{1}{2}
&2b^2-2b-\frac{1}{2}\\
\sqrt{b-b^2}&(2b-1)\sqrt
{b-b^2}\end{pmatrix}
\begin{pmatrix}u\\v
\end{pmatrix}
=\begin{pmatrix}1&bH&b^2-
\frac{1}{2}b\\
0&\sqrt{b-b^2}H&b\sqrt{b-b^2}
\end{pmatrix}.
$$
Hence ${\sigma}^{b}$ also 
satisfies the condition~2 
of Proposition~\ref{geometric1} 
and ${\sigma}^{b}\in\Stab(\PP^2)$ 
is geometric. The proof of 
Proposition~\ref{geometric1}
implies that there exists a 
lift $g\in\grp$ of 
$T\in\GL^+(2,\R)$ such that 
$\pi({\sigma}^{b} g)=
\pi({\sigma}^{b})T$ and  
\begin{equation}\label{ast}
{\sigma}^{b} g=
\sigma_{(bH,tH)},
\end{equation} 
where we put $t=\sqrt{b-b^2}$.
We fix $\alpha\in K(\PP^2)$ with 
$\ch(\alpha)=(r,sH,\ch_2)$, $0<s\le r$.
By the remark after Main~Theorem~\ref{our} 
we may assume that $\ch_2\le\frac{1}{2}$.
We choose $0<b<\frac{s}{r}$ such that 
$\alpha\in K(\PP^2)$ and $\sigma_{(bH,tH)}=
(Z_{(bH,tH)},\A_{(bH,tH)})$ satisfy the conditions in 
Theorem~\ref{limit};
\begin{equation}
   \begin{split}
        0<\e=\Im Z_{(bH,tH)}(\alpha)&=s-rb\le
        \min\left\{t=\sqrt{b-b^2},\ \frac{1}{r}\right\}\\
        \text{and }\hspace{15mm}\Re Z_{(bH,tH)}(\alpha)&=
        -\ch_2+r/2(b-2b^2)+sb\ge 0.\label{sm6}
   \end{split}
\end{equation}
In the following we assume that $s/r-b>0$ is small enough
such that these inequalities 
are satisfied. Then by Corollary~\ref{cor} we have
\begin{equation}\label{isom1}
\M_{\D(\PP^2)}
(\ch(\alpha),\sigma_{(bH,tH)})
\cong\M_{\PP^2}(\ch(\alpha),H).
\end{equation} 

Since $\sigma^bg=\sigma_{(bH,tH)}$,
by (\ref{gisom}) we see that 
the shift functor $\ \cdot \ [n]$ 
gives an isomorphism
\begin{equation}\label{isom3}
\M_{\D(\PP^2)}(\ch(\alpha),\sigma_{(bH,tH)})\cong
\M_{\D(\PP^2)}((-1)^n\ch(\alpha),{\sigma}^{b})\colon
E\mapsto E[n]
\end{equation}
for some $n\in\Z$. 
We show that $n=1$.
First notice that $\alpha=
a_0[F_0]+a_1[F_1]
+a_2[F_2]\in K(\PP^2)$, where
$(a_0,a_1,a_2)\in\Z^3$
is defined by 
\begin{equation*}
\begin{split}
a_0&:=r-\frac{3}{2}s+\ch_2\\
a_1&:=-2s+2\ch_2\\
a_2&:=-\frac{s}{2}+\ch_2.
\end{split}
\end{equation*}
For every $\C$-valued point 
$E$ of $\M_{\D(\PP^2)}
(\ch(\alpha),\sigma_{(bH,tH)})$,
by Lemma~\ref{bei1}~(2), 
$E[n]$ is written as
\begin{equation}
\label{coef}
E[n]\cong\left(\mo_{\PP^2}
^{(-1)^na_0}\to\mo_{\PP^2}
(1)^{(-1)^na_1}\to\mo_{
\PP^2}(2)^{(-1)^na_2}\right)
\in\A_1,
\end{equation}
where $\mo_{\PP^2}(2)^{
(-1)^na_2}$
lies on degree $0$.
The conditions that $0<s\le r$ 
and $\ch_2\le \frac{1}{2}$ imply 
that $a_2\le 0$ and that 
$a_2=0$ if and only if $\ch(
\alpha)=(1,1,\frac{1}{2})$.
In the case $a_2<0$, 
the form~(\ref{coef})
of $E[n]$ implies $n=1$ since $E$ is a 
sheaf. In the case $a_2=0$,
we have $M_{\PP^2}
(\ch(\alpha),H)=\{\mo_{\PP^2}(1)\}$. 
Since $\mo_{\PP^2}(1)[1]\in\A_1$, 
we also have $n=1$.

On the other hand we define 
$\theta^{\alpha}_{Z^{b}}
\colon K(B)\to\R$ by (\ref{det}) using $\varphi_1
\colon K(\PP^2)\cong K(B)$.
Then by Proposition~\ref{algebraic} 
the moduli functor 
$\M_{\D(\PP^2)}(-\ch(\alpha),{\sigma}^{b})$ is 
corepresented by the moduli scheme 
$M_{B}(-\varphi_1(\alpha),\theta^{\alpha}_{Z^{b}})$.
Combining this with the above 
isomorphisms (\ref{isom1})
and (\ref{isom3}) with $n=1$ 
we have an isomorphism
\begin{equation}\label{isom4}
M_{\PP^2}(\ch(\alpha),H)\cong M_{B}(-\varphi_1(\alpha),
\theta^{\alpha}_{Z^{b}})\colon E\mapsto \Phi_1(E[1]).
\end{equation}
Isomorphisms (\ref{isom1}) and 
(\ref{isom3}) hold for moduli functors
corresponding to stable objects.
Hence the isomorphism (\ref{isom4})
keeps open subsets of stable objects.

Finally we see that if 
$s/r-b>0$ is small enough,
this $\theta^{\alpha}_{Z^{b}}$ belongs to 
$C^{\PP^2}_{\varphi_1(\alpha)}$ in the Main~Theorem as
follows. The above 
isomorphism~(\ref{isom4})
implies that if $s/r-b>0$ is 
small enough, $\theta^{\alpha}_{Z^{b}}$
belongs to the same chamber 
$C_{\varphi_1(\alpha)}$. This chamber $C_{\varphi_1(\alpha)}$
satisfies the desired conditions. 
In fact we have 
$\theta^{\alpha}_{Z^{b}}(\varphi_1(\mo_x))> 0$ 
for $b<s/r$ and $\theta^{\alpha}
_{Z^{\frac{s}{r}}}(\varphi_1(\mo_x))=0$, 
furthermore $M_{\PP^2}(
\ch(\alpha),H)\neq\emptyset$ 
implies $M_{B}(-\varphi_1(\alpha),
\theta_1)\neq\emptyset$ 
for $\theta_1\in C_{\varphi_1(\alpha)}$ 
because of the isomorphism 
(\ref{isom4}).    
This completes the proof of 
Main~Theorem~\ref{our}.
\label{prmain}

\subsection{Comparison 
with Le Potier's result}
\label{comp}
In the sequel we show that our 
Theorem~\ref{our} 
implies Main Theorem~\ref{inmain}~(ii), (iii),
in particular, Le Potier's result.
In addition to $\mathfrak{E}_1$,
we consider 
the following full strong 
exceptional collections
on $\PP^2$ 
$$\mathfrak{E}'_1=\left(\mo_{\PP^2}(1),
\mo_{\PP^2}(2),\mo_{\PP^2}(3)\right),
\hspace{15mm}
\mathfrak{E}_0=\left(\mo_{\PP^2}(1),
\Omega^1_{\PP^2}(3),\mo_{\PP^2}(2)
\right),$$
the equivalences $\Phi'_1(\ \cdot\ )=\mathbf R
\Hom_{\PP^2}(\E'_1,\ \cdot\ )$,
$\Phi_0(\ \cdot\ )=\mathbf R
\Hom_{\PP^2}(\E_0,\ \cdot\ )$
between $\D(\PP^2)$and 
$\D(B')$, $\D(B)$
and the induced isomorphisms
$\varphi'_1\colon K(\PP^2)\cong K(B')$, 
$\varphi_0\colon K(\PP^2)\cong K(B)$,
where $\E'_1=\mo_{\PP^2}(1)\oplus
\mo_{\PP^2}(2)\oplus\mo_{\PP^2}(3)$, 
$\E_0=\mo_{\PP^2}(1)\oplus
\Omega^1_{\PP^2}(3)\oplus\mo_{\PP^2}(2)$
and $B'=\End_{\PP^2}(\E_1')$,
$B=\End_{\PP^2}(\E_0)$.
We also recall from \S~\ref{fsecp2} that 
\begin{equation}\label{gen}
\A'_1=\langle \mo_{\PP^2}[2],
\Omega^1_{\PP2}(2)[1],
\mo_{\PP^2}(1)\rangle,
\hspace{5mm}\A_0=\langle
\mo_{\PP^2}(-1)[2],
\mo_{\PP2}[1],\mo_{\PP^2}(1)
\rangle.
\end{equation}
We remark that $\A'_1$ is the left tilt of $\A_1=
\langle\mo_{\PP^2}[2],\mo_{\PP^2}(1)[1],
\mo_{\PP^2}(2)\rangle$ at  
$\mo_{\PP^2}(1)[1]$ and $\A_0$ is the 
left tilt of $\A'_1$ at $\mo_{\PP^2}[2]$.
See \cite{B4} for this terminology and relationship
between tilting and exceptional collections
although we do not use this fact. 

For $\theta\in\Hom_{\Z}(K(\PP^2),\R)$, 
we put $\theta_k:=
\theta\circ{\varphi_k}^{-1}\in
\Hom_{\Z}(K(B),\R)$ for $k=0,1$ 
and $\theta'_1:=\theta\circ
{\varphi'_1}^{-1}\in
\Hom_{\Z}(K(B'),\R)$.
We put 
\begin{equation}\label{vtheta}
\begin{split}
(\theta_k^0,\theta_k^1,\theta_k^2)
&:=(\theta_k(\C v_0),\theta_k(\C v_1),\theta_k(\C v_2))
\hspace{7mm}\text{for }k=0,1,\\  
({\theta'_1}^0,{\theta'_1}^1,{\theta'_1}^2)
&:=(\theta'_1(\C v_0),\theta'_1(\C v_1),\theta'_1(\C v_2)).
\end{split}
\end{equation}
For any $B$-module $N$ and $B'$-module
$M$, we have  
\begin{equation*}
\begin{split}
\theta_k(N)&=\theta_k^0\dim_{\C}(Nv^{\ast}_0)
+\theta_k^1\dim_{\C}(Nv^{\ast}_1)+\theta_k^2\dim
_{\C}(Nv^{\ast}_2)\hspace{7mm}\text{for }k=0,1,\\
\theta'_1(M)&={\theta'_1}^0\dim_{\C}(Mv^{\ast}_0)
+{\theta'_1}^1\dim_{\C}(Mv^{\ast}_1)+
{\theta'_1}^2\dim_{\C}(Mv^{\ast}_2).
\end{split}
\end{equation*}
By abbreviation we denote this by 
$\theta_k=(\theta_k^0,\theta_k^1,
\theta_k^2)$ and $\theta'_1=({\theta'_1}^0,
{\theta'_1}^1,{\theta'_1}^2)$. 
It is also convenient to write the 
following equality 
\begin{equation}\label{simp}
\begin{split}
(\theta_k^0,\theta_k^1,
\theta_k^2)&=\left(\theta(\mo_{\PP^2}(k-1)[2]),
\theta(\mo_{\PP^2}(k)[1]\right),
\theta(\mo_{\PP^2}(k+1)))\hspace{7mm}\text{
for }k=0,1,\\
({\theta'_1}^0,{\theta'_1}^1,
{\theta'_1}^2)&=\left(\theta(\mo_{\PP^2}[2]),
\theta(\Omega_{\PP^2}(2)[1]),
\theta(\mo_{\PP^2}(1))\right).
\end{split}
\end{equation}

\begin{prop}\label{tilt}
Let $\theta\colon K(\PP^2)\to\R$ 
be an additive function with $\theta_1=
(\theta_1^0,\theta_1^1,\theta_1^2)$ and 
$\alpha\in K(\PP^2)$ with 
$\theta(\alpha)=0$. 
If $\theta_1^0,\theta_1^1<0$, 
then equivalences $\Phi'_1\circ\Phi_1^{-1}
\colon\D({B})\cong\D(B')$ and 
$\Phi_0\circ{\Phi'_1}^{-1}\colon
\D(B')\cong\D(B)$ between derived 
categories induce the isomorphisms
$$M_{B}(\varphi_1(\alpha),\theta_1)
\cong M_{B'}(\varphi'_1(\alpha),\theta'_1)\cong 
M_{B}(\varphi_0(\alpha),\theta_0).$$
These isomorphisms keep
open subsets of stable modules.
\end{prop}

We only show the first isomorphism using the
assumption that $\theta_1^1<0$. The
other assumption that $\theta_1^0<0$
is used for the second isomorphism. \\

\noindent \emph{\bf{Step 1.}}\ 
The assumption 
$\theta_1^1<0$ implies that 
$\Phi'_1\circ\Phi_1^{-1}(N)\in$ mod-${B'}$
for any $N\in M_{B}(\varphi_1(\alpha),\theta_1)$.
\proof
We take 
$E\in\A_1$ such that $\Phi_1(E)=N$. Then
the decomposition of $N=\mathbf R\Hom_{\PP^2}
(\E_1,E)$ is given by 
\begin{equation}\label{module}
\begin{split}
Nv^{\ast}_0&=\mathbf R\Hom_{\PP^2}(\mo_{\PP^2}(3),E)\\
Nv^{\ast}_1&=\mathbf R\Hom_{\PP^2}(\Omega^1_{\PP^2}(4),E)\\
Nv^{\ast}_2&=\mathbf R\Hom_{\PP^2}(\mo_{\PP^2}(2),E),
\end{split}
\end{equation}
and $\gamma_i^\ast|_N=p_i^{\ast}$, 
$\delta^{\ast}_j|_N=q_j^{\ast}$ from 
(\ref{final}).
On the other hand, we have 
\begin{equation}\label{deco}
\begin{split}
\Phi'_1\circ\Phi^{-1}_1(N)&=\mathbf 
R\Hom_{\PP^2}(\mathcal{E}'_1,E)\\
&=\mathbf R\Hom_{\PP^2}(\mo_{\PP^2}(3),E)\oplus
\mathbf R\Hom_{\PP^2}(\mo_{\PP^2}(2),E)\oplus
\mathbf R\Hom_{\PP^2}(\mo_{\PP^2}(1),E).
\end{split}
\end{equation}
The fact that $N\in$ mod-$B$ and (\ref{module})
implies 
$$\mathbf R^i\Hom_{\PP^2}(\mo_{\PP^2}(3),E)=
\mathbf R^i\Hom_{\PP^2}(\mo_{\PP^2}(2),E)=0$$
for $i\neq 0$.
From the exact sequence  
\begin{equation}\label{kos}
0\to\mo_{\PP^2}(1)\xrightarrow{\Sigma
z_i\otimes e_i}
\mo_{\PP^2}(2)\otimes V
\xrightarrow{q_i\otimes e_i^\ast}
\Omega^1_{\PP^2}(4)\to 0,
\end{equation}
we have an isomorphism of 
complexes in $\D(\PP^2)$
\begin{equation}\label{proj}
\mo_{\PP^2}(1)\cong\left(
\mo_{\PP^2}(2)\otimes V \xrightarrow{
\Sigma q_i\otimes e^{\ast}_i} 
\Omega^1_{\PP^2}(4)\right),
\end{equation}
where $\mo_{\PP^2}(2)\otimes V$ lies 
on degree $0$.
By applying Lemma~\ref{bei1}~(1) 
to $(\ref{proj})$ and $E\in\A_1$,
we have an isomorphism in $\D(\C)$
\begin{equation}\label{mid}
\mathbf R\Hom_{\PP^2}(\mo_{\PP^2}(1),E)
\cong\left(Nv^{\ast}_1\xrightarrow{\delta^{\ast}_V}
(Nv^{\ast}_2)\otimes V\right),
\end{equation}
where $(Nv^{\ast}_2)\otimes V$ lies on 
degree $0$ and $\delta^{\ast}_V=
\delta^{\ast}_0\otimes e_0+
\delta^{\ast}_1\otimes e_1+
\delta^{\ast}_2\otimes e_2$.
Hence  
$\Phi'_1\circ\Phi_1^{-1}(N)$ belongs
to $\text{mod-}{B'}$ if and only if 
$$\ker \delta^{\ast}_V=\mathbf R^{-1}\Hom_{\PP^2}
(\mo_{\PP^2}(1),E)=0.$$
However if $\ker \delta^{\ast}_V\neq 0$, we
can view $\ker \delta^\ast_V$ as a
submodule $N'$ of $N$ with 
$N'v^{\ast}_0=N'v^{\ast}_2=0$ and $N'v^{\ast}_1=\ker \delta^\ast_V$.
This contradicts $\theta_1$-semistability 
of $N$ since $\theta_1(\ker \delta^{\ast}_V)
=\theta_1^1\cdot\dim_{\C} (\ker \delta^{\ast}_V)<0$.  
\endproof

\noindent \emph{\bf{Step 2.}} 
For any $N\in M_{B}(\varphi_1(\alpha),\theta_1)$,
$\theta_1$-(semi)stability of $N$
implies $\theta'_1$-(semi)stability 
of $M:=\Phi'_1\circ\Phi_1^{-1}(N)
\in\text{mod-}{B'}$.
\proof
We recall that $v_i\in \C Q/J'$ correspond to $\id_{\mo_{\PP^2}(3-i)}
\in B'$
for $i=0,1,2$ via the isomorphism (\ref{final}). 
Hence by (\ref{module}), (\ref{deco}) 
and (\ref{mid}) we have 
\begin{equation}\label{decom}
Mv^{\ast}_0=Nv^{\ast}_0, \ \ Mv^{\ast}_1=Nv^{\ast}_2,\ \ 
Mv^{\ast}_2=\coker 
\delta^{\ast}_V.
\end{equation}
Since $z_i=p_{i+2}\circ q_{i+1}
\in\Hom_{\PP^2}(\mo_{\PP^2}(2),
\mo_{\PP^2}(3))$, $\gamma^{\ast}_i|_M\colon
Mv^{\ast}_0\to Mv^{\ast}_1$
is defined by 
\begin{equation*}
\gamma^{\ast}_i|_M:=\delta^{\ast}_{i+1}|_N
\circ\gamma^{\ast}_{i+2}|_N\colon Nv^{\ast}_0\to 
Nv^{\ast}_2.
\end{equation*}
Via the isomorphism (\ref{proj}), 
homomorphisms $z_i
\colon\mo_{\PP^2}(1)\to\mo_{\PP^2}
(2)$ correspond to homotopy 
classes of homomorphisms $\id_{\mo_{\PP^2}(2)}
\otimes e_i^{\ast}
\colon \mo_{\PP^2}(2)\otimes V\to\mo_{\PP
^2}(2)$ in 
$$\Hom_{\D(\PP^2)}(\mo_{\PP^2}(1),
\mo_{\PP^2}(2))\cong \coker\left(\Hom_{\PP^2}
(\Omega^1_
{\PP^2}(4), \mo_{\PP^2}(2))\to
\Hom_{\PP^2}
(\mo_{\PP^2}(2)\otimes V, \mo_{\PP^2}(2)\right)$$
for $i=0,1,2$. Hence $\delta^{\ast}_j|_M
\colon Mv^{\ast}_1\to Mv^{\ast}_2$
is defined by 
\begin{equation*}
\delta^{\ast}_j|_M\colon Nv^{\ast}_2\xrightarrow
{\id_{Nv^{\ast}_2}\otimes e_j} (Nv^{\ast}_2)\otimes V
\to\coker 
\delta^{\ast}_V,
\end{equation*}
where $(Nv^{\ast}_2)\otimes V\to\coker 
\delta^{\ast}_V$ is a natural 
surjection.

Conversely from this description 
we see easily that 
the above $B$-module $N$ is 
reconstructed from the ${B'}$-module 
$M=\Phi'_1\circ\Phi_1^{-1}(N)$ as follows. 
We define 
\begin{equation}\label{delta}
{\delta^{\ast}}^V:=
\Sigma_i(\delta_i^{\ast}|_M)\otimes e^{\ast}_i
\colon (Mv^{\ast}_1)\otimes V\to Mv^{\ast}_2.
\end{equation}
We put 
\begin{equation}
\label{decon}
Nv^{\ast}_0:=Mv^{\ast}_0,\ \ Nv^{\ast}_1:=\ker {\delta^{\ast}}^V,
\ \ Nv^{\ast}_2:=Mv^{\ast}_1
\end{equation}
and define $\gamma^{\ast}_i|_N\colon Nv^{\ast}_0\to Nv^{\ast}_1$ and 
$\delta^{\ast}_j|_N\colon Nv^{\ast}_1\to Nv^{\ast}_2$ by
\begin{equation}\label{action}
\begin{split}
\gamma^{\ast}_i|_N&:=
(\gamma_{i+1}^{\ast}|_M)\otimes 
e_{i+2}-(\gamma^{\ast}_{i+2}
|_M)\otimes e_{i+1}\colon Mv^{\ast}_0
\to\ker{\delta^{\ast}}^V,\\
\delta^{\ast}_j|_N&\colon
\ker{\delta^{\ast}}^V\subset 
(Mv^{\ast}_1)\otimes V\xrightarrow{\id_
{Mv^{\ast}_1}\otimes e_j^{\ast}}Mv^{\ast}_1.
\end{split}
\end{equation} 
Imitating this, for any 
${B'}$-submodule 
$M'$ of $M$ we construct 
an ${B}$-submodule 
$N'$ of $N$ by (\ref{delta}), 
(\ref{decon}) and (\ref{action}) 
with $Mv^{\ast}_i$ and $Nv^{\ast}_j$ replaced
by $M'v^{\ast}_i$ and $N'v^{\ast}_j$. However 
in this case 
$${\delta^\ast}^V
\colon (M'v^{\ast}_1)\otimes V\to M'v^{\ast}_2$$ 
is not necessarily
surjective. Hence we have 
$$\dim_{\C}(N'v^{\ast}_1)=\dim_{\C} 
\ker\left({\delta^\ast}^V|
_{(M'v^{\ast}_1)\otimes V}\right)
\ge 3\dim_\C (M'v^{\ast}_1)-\dim_\C (M'v^{\ast}_2).$$
Hence the assumption that 
$\theta_1^1<0$ and the following
equality by (\ref{simp})
$$(\theta_1^0,
\theta_1^1,\theta_1^2)
\begin{pmatrix}
1&0&0\\
0&3&-1\\
0&1&0
\end{pmatrix}=
({\theta'_1}^0,
{\theta'_1}^1,{\theta'_1}^2)$$ 
implies $\theta_1(N')\le\theta'_1(M')$.
Thus $\theta_1$-(semi)stability of
$N$ implies $\theta'_1$-(semi)stability
of $M$ and we have
$$\Phi'_1\circ\Phi_1^{-1}(M_{B}(\varphi_1(\alpha),\theta_1))
\subset M_{B'}(\varphi'_1(\alpha),\theta'_1).$$
The proof of the opposite inclusion is similar 
and we leave it to the readers.
\endproof

If we assume $\ch_2<\frac{1}{2}$, 
the chamber $C^{\PP^2}_{\varphi_1(\alpha)}
\subset\varphi_1(\alpha)^{\perp}$ 
defined in Section~\ref{our} 
intersect with the region 
defined by the inequalities 
$\theta_1^0,\theta_1^1<0$. 
Hence from the above proposition
and Theorem~\ref{our} we
have isomorphisms
\begin{equation}\label{potier}
M_{\PP^2}(\ch(\alpha),H)\cong 
M_{B'}(-\varphi'_1(\alpha),\theta'_1)
\colon E\mapsto \Phi'_1(E[1])
\end{equation} 
\begin{equation}
M_{\PP^2}(\ch(\alpha),H)\cong 
M_{B}(-\varphi_0(\alpha),\theta_0)
\colon E\mapsto \Phi_0(E[1])
\end{equation} 
for $\alpha\in K(\PP^2)$ with
$0<c_1(\alpha)\le \rk(\alpha)$,
$\ch_2<\frac{1}{2}$ and 
$\theta\colon K(\PP^2)\to\R$
satisfying $\theta_1\in C^{\PP^2}_{\varphi_1(\alpha)}$ with 
$\theta^0_1,\theta^1_1<0$.
This completes the proof of Main Theorem~\ref{inmain}.
(\ref{potier}) was obtained by 
Le Potier \cite{P}.

\section{Computations of the wall-crossing}
\label{cwc}

In this section, we identify 
the Hilbert schemes of
points on $\PP^2$ 
$$(\PP^2)^{[n]}:=\{\mathcal{I}
\subset\mo_{\PP^2}\mid
\text{Length}(\mo_{\PP^2}/
\mathcal{I})=n\}$$
with the moduli spaces
$M_{B}(-\varphi_0(\alpha),\theta_0)\cong
M_{B}(-\varphi_1(\alpha),\theta_1)$ 
by Theorem~\ref{our} and 
Proposition~\ref{tilt}, where 
$\alpha\in K(\PP^2)$ with $\ch(\alpha)
=(1,1,\frac{1}{2}-n)$, 
$\theta_1\in C^{\PP^2}_{\varphi_1(\alpha)}$ and 
$\theta_0=\theta_1\circ\varphi_1\circ
\varphi_0^{-1}$. 
We study the wall-crossing phenomena 
of the Hilbert schemes of
points on $\PP^2$ via this identification.

\subsection{Geometry of 
Hilbert schemes of points 
on $\PP^2$}
\label{eff}
We recall the geometry of Hilbert schemes 
of points on $\PP^2$ (cf.~\cite{LQZ}). 
Let $\ell$ be a line in $\PP^2$, and 
$x_1,\dots ,x_{n-1}\in \PP^2$ be 
distinct fixed points in $\ell$. Let
$$M_2(x_1)=\big\{\,\xi\in (\PP^2)^{[2]}
\mid\text{Supp}(\xi)={x_1}\big\}$$ 
be the punctual Hilbert 
scheme parameterizing length-$2$ 
$0$-dimensional subschemes 
supported at $x_1$. 
It is known that $M_2(x_1)\cong \PP^1$. 
Let $N_1((\PP^2)^{[n]})$ be the 
$\R$-vector space of numerical equivalence 
classes of one-cycles on $(\PP^2)^{[n]}$.
We define two curves $\beta_n$ and 
$\zeta_{\ell}$ in $(\PP^2)^{[n]}$
as elements in $N_1((\PP^2)^{[n]})$
by the following formula
\begin{equation}
\begin{split}\label{defn}
\beta_n&:=\big\{\,\xi+x_2+\dots 
+x_{n-1}\in (\PP^2)^{[n]}\mid
\xi\in M_2(x_1)\big\}\\
\zeta_\ell&:=\big\{\,x+x_1+\dots 
+x_{n-1}\in (\PP^2)^{[n]}
\mid x\in \ell\big\}.
\end{split}
\end{equation}
The definition
of $\beta_n$ and $\zeta_\ell$ 
does not depend on the choice of 
a line $\ell$ on $\PP^2$
and points $x_1,\ldots,x_{n-1}$
on $\ell$ (cf. \cite[Theorem~3.2 and
Theorem~5.1]
{LQZ}).
We define a cone $\NE((\PP^2)^{[n]})$ in 
$N_1((\PP^2)^{[n]})$ by
$$\NE((\PP^2)^{[n]}):=\left\{\Sigma a_i[C_i]
\mid C_i\subset (\PP^2)^{[n]}
\text{ an irreducible curve, }
a_i\ge 0\right\}$$
and $\overline{\NE}((\PP^2)^{[n]})$ to be 
its closure. 
\begin{thm}\label{LQZ}\emph{\bf{\cite[Theorem~4.1]{LQZ}}}
$\overline{\NE}((\PP^2)^{[n]})$ 
is spanned by 
$\beta_n$ and $\zeta_\ell$.
\end{thm}

Let $S^n(\PP^2)$ be the $n$th symmetric 
product of $\PP^2$, that is, $S^n(\PP^2):=(\PP^2)^n/
\mathfrak{S}_n$, where $\mathfrak{S}_n$ is the 
symmetric group of degree $n$. The Hilbert-Chow 
morphism $\pi\colon (\PP^2)^{[n]}\to S^n(\PP^2)$ is 
defined by $\pi(\mathcal{I})=\text{Supp}(\mo_{\PP^2}/
\mathcal{I})\in S^n(\PP^2)$ for every $\mathcal{I}\in 
(\PP^2)^{[n]}$. The morphism $\pi$ is the contraction of 
the extremal ray $\R_{>0}\beta_n$. 

Denote by $\psi\colon (\PP^2)^{[n]}
\to Z$ the contraction 
morphism of the extremal ray 
$\R_{>0}\zeta_\ell$. 
In the case $n=2$, 
$\psi\colon (\PP^2)^{[2]}\to Z$ 
coincide with the morphism
$\text{Hilb}^2(\PP((T_{(\PP^2)^{\ast}})
^{\ast}))\to(\PP^2)^{\ast}$ up
to isomorphism, where 
$\text{Hilb}^2(\PP((T_{(\PP^2)^{\ast}})
^{\ast}))$ is the relative Hilbert scheme. 
In the case $n=3$, 
$\psi\colon (\PP^2)^{[3]}\to Z$ is 
a divisorial contraction. 
In the case $n\ge 4$, 
$\psi\colon (\PP^2)^{[n]}\to Z$ 
is a flipping 
contraction.

\subsection{Wall-Crossing of the 
Hilbert schemes of points on $\PP^2$}
\label{wallc}
We take $\alpha\in K(\PP^2)$ with 
$\ch(\alpha)=(r,1,\frac{1}{2}-n)$ and
assume that $n\ge 1$. By (\ref{dim}),
we have $\dimv(-\varphi_1(\alpha))=(n-r+1,2n+1,n)$.
For $b\in\R$ with 
$0<b<\frac{1}{r}$ we put $t=\sqrt{
b-b^2}$. 
From (\ref{isom3}) and 
Proposition~\ref{algebraic}, 
we have isomorphisms
\begin{equation}
\label{ident1}
{}^{sh}\M_{\D(\PP^2)}(\ch(\alpha),\sigma_{(bH,
tH)})
\cong
{}^{sh}\M_{\D(\PP^2)}(-\ch(\alpha),\sigma^b)
\colon E\mapsto E[1]
\end{equation} 
\begin{equation}\label{ident2}
{}^{sh}\M_{\D(\PP^2)}(-\ch(\alpha),\sigma^b)
\cong {}^{sh}\M_{B}(-\varphi_1(\alpha),\theta^{\alpha}_{Z^b})
\colon E[1]\mapsto \Phi_1(E[1]),
\end{equation}
where $\sigma^b$ is
defined by (\ref{def1}) and $\theta^{\alpha}_{Z^b}$
is defined by (\ref{det}) using $\varphi_1\colon
K(\PP^2)\cong K(B)$. 
We recall that from \S~\ref{proof}, if 
$\frac{1}{r}-b_0>0$ is small enough,
then $M_{\PP^2}(\ch(\alpha),H)$ corepresents 
${}^{sh}\M_{\D(\PP^2)}(\ch(\alpha),\sigma_{(b_0H,
t_0H)})$, where $t_0:=\sqrt{b_0-b_0^2}$. 
We have $\theta^{\alpha}_{Z^{b_0}}\in C^{\PP^2}_{\varphi_1(\alpha)}$ 
and the isomorphism 
$$M_{\PP^2}(\ch(\alpha),H)\cong M_{B}(-
\varphi_1(\alpha),\theta^{\alpha}_{Z^{b_0}})$$ in
Theorem~\ref{our}.
In fact the following lemma holds.
\begin{lem}\label{chamber}
We have $\R_{>0}\theta^\alpha_{Z^0}+\R_{>0}\theta
^{\alpha}_{Z^\frac{1}{r}}\subset C^{\PP^2}_{\varphi_1(\alpha)}$,
that is, the moduli functor 
$\M_{\D(\PP^2)}
(\ch(\alpha),\sigma_{(bH,tH)})$
does not change as $b$ moves in 
the interval $(0,\frac{1}{r})$. 
\end{lem}
\proof
We assume that there exists a $\C$-valued 
point $E$ of $\M_{\D(\PP^2)}(\ch(\alpha),\sigma_{(b_0H,
t_0H)})$ 
such that $E$ is not $\sigma_{(b_1H,
t_1H)}$-semistable 
for some $b_1\in(0,\frac{1}{r})$, where 
we put $t_1:=\sqrt{b_1-b_1^2}$.
From (\ref{ident1}) and (\ref{ident2}), 
$\sigma_{(bH,tH)}$-semistability for $E$
and $\theta^{\alpha}_{Z^b}$-semistability for 
$\Phi_1(E[1])$ are equivalent for $b\in
(0,\frac{1}{r})$.
Using the notation (\ref{vtheta}) in
\S~\ref{comp}, $\theta^{\alpha}_{Z^b}$ is computed from
(\ref{def1}) and (\ref{simp}) as follows:
$$\theta^{\alpha}_{Z^b}
=(1-b)(0,-n,2n+1)+b(-n,0,n+1-r)\in
\Hom_{\Z}(K(B),\R)\cong\R^3.$$
If we fix any $\beta\in K(B)$, then 
$\theta^{\alpha}_{Z^b}(\beta)$ is a monotonic
function for $b$.
Hence we may assume that such a 
real number $b_1$ is small enough.

We take the $\sigma_{(b_1H,t_1H
)}$-semistable factor 
$G$ of $E$ 
with the smallest slope $\mu_
{\sigma_{(b_1H,t_1H)}}(G)$ and the 
exact sequence in $\A_{(b_1H,t_1H)}$
\begin{equation}\label{aex}0\to F
\to E \to G\to 0,
\end{equation}
where $F$ is a nonzero object 
of $\A_{(b_1H,t_1H)}$.
From (\ref{aex}) we see that $F$ 
is a sheaf since $E$ is a sheaf
and $\H^i(G)=0$ for $i\neq 0,-1$.
From (\ref{ident1}) we have 
$E[1]\in\A_1$. By the uniqueness
of Harder-Narasimhan filtration we see
that $G[1]$ and $F[1]$ also belong to $\A_1$.
Hence from the exact sequence (\ref{aex}), we
see that dimension vectors of $B$-modules
$\Phi_1(G[1])$ and $\Phi_1(F[1])$ are
bounded from above by $\dimv(-\varphi_1(\alpha))$. 
In particular there exists a 
bound of $\rk(F)$ and $\rk(G)$ 
independent of the choice of 
$E$ and $b_1$. The inequality 
$0<\Im Z_{(b_1H,t_1H)}(F)=
t_1(c_1(F)-r(F)b_1)
<\Im Z_{(b_1H,t_1H)}(E)$ implies that 
$0<c_1(F)\le c_1(E)=1$ since we can take 
arbitrary small $b_1>0$ and $\rk(F)$ is
bounded from above. So we have $c_1(F)=1$ and 
$c_1(G)=c_1(E)-c_1(F)=0$.

We put $I:=\im(F\to E)$.
Since $F\to I$ is surjective we have $0<\mu_{H\text{-min}}(F)
\le\mu(I)$. Furthermore since 
$E$ is Gieseker-semistable,
we have $\mu(I)\le \mu(E)=
\frac{1}{r}$. Hence 
$\rk(I)=r$, $c_1(I)=1$
and $\H^0(G)$ is a $0$-dimensional
sheaf. Since $G[1]\in\A_1$, by 
Lemma~\ref{bei1}~(2)
we have an isomorphism
$$G[1]\cong\left(\mo_{\PP^2}
^{\oplus a_0}\to\mo_{\PP^2}
(1)^{\oplus a_1}\to\mo_{
\PP^2}(2)^{\oplus a_2}
\right),$$
where $( a_0, a_1, a_2)=-
r(G)(1,0,0)-\ch_2(G)(1,2,1)\in
\Z_{\ge0}^3$. Hence $\ch_2(G)$ 
must be non-positive and 
$\ch_2(G)=0$ if and only if $G[1]\cong\mo
_{\PP^2}^{\oplus a_0}[2]$. 
In this case, we have $\theta^{\alpha}_{Z^{b_1}}
(\Phi_1(G[1]))=-nb_1 a_0<0$ and
$\Phi_1(G[1])$ does not break $\theta^{\alpha}
_{Z^{b_1}}$-semistability of $\Phi_1(E[1])$.
This contradicts the choice of $G$. 
We have $\ch_2(\H^{-1}(G))=-\ch_2(G)+
\ch_2(\mathcal{H}^0(G))>0$.
On the other hand, we have
$c_1(\H^{-1}(G))=-c_1(G)+c_1(\mathcal{H}^0(G))=0$ 
and from $G\in\A_{(b_1H,t_1H)}$ we have 
$\mu_{H\text{-max}}(\H^{-1}(G))\le 0$ for small enough
$b_1>0$. Hence $\H^{-1}(G)$
is $\mu_H$-semistable and 
satisfy the inequality 
$-2r(\H^{-1}(G))\ch_2(\H^{-1}(G))\ge 0$
by Theorem~\ref{bgine}.
This is a contradiction.
\endproof

In the following we consider 
the case $r=1$. We fix 
$\alpha\in K(\PP^2)$ with 
$\ch(\alpha)=(1,1,\frac{1}{2}-n)$,
$n\ge 1$ and $\theta_1\in C^{\PP^2}_{\varphi_1(\alpha)}$.
Tensoring by $\mo_{\PP^2}(1)
=\mo_{\PP^2}(H)$ does not change 
Gieseker-semistability 
of torsion free sheaves on $\PP^2$ 
and induces an automorphism of 
$K(\PP^2)$ sending $\Hat{\alpha}$ with 
$\ch(\Hat{\alpha})=(1,0,-n)$
to $\alpha$. Since by 
definition $(\PP^2)^{[n]}
=M_{\PP^2}(\ch(\Hat{\alpha}),H)$, 
we have an isomorphism 
$$(\PP^2)^{[n]}\cong M_{\PP^2}
(\ch(\alpha),H)\colon \mathcal{I}
\mapsto\mathcal{I}(1).$$
On the other hand, by 
Theorem~\ref{our} and
Proposition~\ref{tilt},
we have 
isomorphisms
$$\Phi_k(\ \cdot\ [1])\colon
M_{\PP^2}(\ch(\alpha),H)\cong
M_{B}(-\varphi_k(\alpha),\theta_k)$$
for $k=0,1$, where $\theta_0=\theta_1
\circ\varphi_1\circ\varphi_0^{-1}$.
In what follows,
we often use these 
identifications
$$(\PP^2)^{[n]}\cong M_{B}(
-\varphi_k(\alpha),\theta_k)\colon
\mathcal{I}\mapsto\Phi_k(\mathcal{
I}(1)[1]),\hspace{.5cm} \text{ and } 
\hspace{.5cm}
\Phi_k\colon\A_k\cong \text{mod-}B.$$
For any $0$-dimensional subscheme
$Z$ of $\PP^2$, $\mathcal{I}_Z$ 
denotes the ideal of $Z$, that is, 
the structure sheaf $\mo_Z$
is defined by $\mo_Z:=\mo_{\PP^2}/
\mathcal{I}_Z$. If the length of $Z$
is $n$, then $\mathcal{I}_Z$ is 
an element of $(\PP^2)^{[n]}$.

We recall that
\begin{equation}\label{regen}
\A_1=\langle\mo_{\PP^2}[2],
\mo_{\PP^2}(1)[1],\mo_{\PP^2}(2)\rangle,
\ \ \ \ \A_0=\langle\mo_{\PP^2}(-1)[2],
\mo_{\PP^2}[1],\mo_{\PP^2}(1)\rangle,
\end{equation} 
$$\dimv(-\varphi_1(\alpha))=(n,2n+1,n),
\ \ \ \ \dimv(-\varphi_0(\alpha))=
(n,2n,n-1).$$
For $b\in\R$, we put
\begin{equation}\label{theta1}
\theta(b)_1:=(1-b)(0,-n,2n+1)+b(-n,0,n)
\in\Hom_{\Z}(K(B),\R)
\end{equation}
\begin{equation}\label{theta2}
\theta(b)_0:=(1-b)(-n+1,0,n)+b(-2n,n,0)
\in\Hom_{\Z}(K(B),\R).
\end{equation}
If $0<b<1$, by (\ref{def1}) and (\ref{simp}) 
we have $\theta(b)_1=\theta^{\alpha}_{Z^b}$
and $\theta(b)_0=\theta^{\alpha}_{Z^b}
\circ\varphi_1\circ\varphi_0^{-1}$.
By Lemma~\ref{chamber}, we have 
$\R_{>0}{\theta(0)_1}
+\R_{>0}{\theta(1)_1}\subset C^{\PP^2}_{\varphi_1(\alpha)}$ 
in $\varphi_1(\alpha)^{\perp}$.
We define a wall-and-chamber structure
on $\varphi_0(\alpha)^{\perp}$ as in \S~\ref{main}
and take the chamber $C^{\PP^2}_{\varphi_0(\alpha)}$ on 
$\varphi_0(\alpha)^{\perp}$ containing $\R_{>0}{\theta(0)}_0
+\R_{>0}{\theta(1)}_0$. 

\begin{lem}\label{C_0}
The following hold.\\
$(1)$ $\R_{>0}\theta(0)_1
+\R_{>0}\theta(1)_1=C^{\PP^2}_{\varphi_1(\alpha)}$
for $n\ge 1$.\\
$(2)$ $\R_{>0}\theta(0)_0
+\R_{>0}\theta(1)_0=C^{\PP^2}_{\varphi_0(\alpha)}$
for $n\ge 2$.
\end{lem}
\proof
It is enough to show that 
$\theta(0)_k$ and 
$\theta(1)_k$ lie on 
walls on $\varphi_k(\alpha)^{\perp}$ for
$k=0,1$.\\
(1) Any $B$-module $N$ with $[N]=\varphi_1(\alpha)$ 
has a surjection
$N\to\C v_0$ and $\theta(0)_1(\C v_0)=0$.
Thus $\theta(0)_1$ lies on a wall on $\varphi_1(\alpha)^
{\perp}$. 
We take any element $\mathcal{I}_Z
\in(\PP^2)^{[n]}$. We have an
exact sequence 
\begin{equation}\label{stex}
0\to\mathcal{I}_Z\to\mo_{\PP^2}\to
\mo_Z\to0.
\end{equation}
$\mo_Z$ can be obtained by extensions of 
$\{\mo_x\mid
x\in\text{Supp}(Z)\}$. Since $\mo_x$ belongs
to $\A_1$ by (\ref{skys}), we have
$\mo_Z\in\A_1$. From (\ref{stex}), 
tensoring by $\mo_{\PP^2}(1)$ 
we have an exact sequence in $\A_1$
$$0\to\mo_Z\to\mathcal{I}_Z(1)[1]\to
\mo_{\PP^2}(1)[1]\to0.$$
Furthermore we have 
$\theta(1)_1(\Phi_1(\mo_Z))=0$, since 
$\dimv(\Phi_1(\mo_{x}))=(1,2,1)$ and 
$\theta(1)_1(\Phi_1(\mo_x))=0$ for any closed 
point $x\in\PP^2$ by (\ref{theta1}).
Thus $\theta(1)_1$ also lies on a wall on 
$\varphi_1(\alpha)^{\perp}$.\\ 
(2) Any $B$-module $N$ with $[N]=\varphi_0(\alpha)$
has a submodule $\C v_2$.
Since $\theta(1)_0(\C v_2)=0$,
$\theta(1)_0$ lies on a wall on 
$\varphi_0(\alpha)^{\perp}$.
On the other hand, 
for any line $\ell$ on $\PP^2$ 
we take an element 
$\mathcal{I}_Z$ of $\zeta_{\ell}
$. 
Since $Z$ is a closed subscheme
of $\ell$ by the definition 
$(\ref{defn})$, we have a diagram:
$$\xymatrix{
0\ar[r]&\mathcal{I}_{Z}\ar[r]
&\mo_{\PP^2}\ \ar[r]
&\mo_Z\ar[r]&0\\
0\ar[r]&\mo_{\PP^2}(-1)\
\ar[u]\ar[r]&\mo_{\PP^2}\ar[r]
\ar@{=}[u]&\mo_{\ell}\ar[u]\ar[r]&0.
}$$
Hence tensoring by $\mo_{\PP^2}(1)$,
we get an
exact sequence in $\Coh(\PP^2)$
$$0\to\mo_{\PP^2}\to\mathcal{I}_Z(1)
\to\mo_{\ell}(-n+1)\to 0,$$ 
where $\mo_{\ell}(-n+1)=\ker\left
(\mo_{\ell}(1)\to\mo_Z\right)$. 
This gives a distinguished 
triangle in $\D(\PP^2)$
\begin{equation}\label{0ex}
\mo_{\PP^2}[1]\to\mathcal{I}_Z(1)[1]
\to\mo_{\ell}(-n+1)[1]\to \mo_{\PP^2}[2].
\end{equation}
We show that this gives an exact sequence in 
$\A_0$. 
It is enough to show that $\mo_{\ell}(-n+1)[1]
\in\A_0$.
An exact sequence in $\Coh(\PP^2)$
$$0\to\mo_{\PP^2}(-1)\to\mo_{\PP^2}
\to\mo_{\ell}\to 0$$
implies that $\mo_{\ell}[1]\in\A_0$ from
(\ref{regen}).
For an integer $m>0$ and a
closed point $x$ in $\ell$, we consider 
an exact sequence in $\Coh(\PP^2)$
$$0\to\mo_{\ell}(-m)\to\mo_{\ell}
(-m+1)\to\mo_{x}\to 0.$$
This gives a distinguished triangle in
$\D(\PP^2)$
$$\mo_{x}\to\mo_{\ell}(-m)[1]
\to\mo_{\ell}(-m+1)[1]\to\mo_{x}[1].$$
Since $\mo_{x}$ belongs to $\A_0$ as in 
Lemma~\ref{geom2}, by induction on $m$
we have $\mo_{\ell}(-m)[1]\in\A_0$ 
for any $m\ge 0$.
Since $\theta(0)_0(\varphi(\mo_{\PP^2}[1]))=0$,
$\mathcal{I}_Z(1)[1]$ and the subobject
$\mo_{\PP^2}[1]$ define a wall 
$\R_{\ge 0}\theta(0)_0$
on $\varphi_0(\alpha)^{\perp}$.
\endproof
 
We take the chamber $C_{\varphi_1(\alpha)}^+
\neq C^{\PP^2}_{\varphi_1(\alpha)}$ in $\varphi_1(\alpha)
^\perp$ sharing the wall $\R_{\ge0}\theta(1)_1$ 
with $C^{\PP^2}_{\varphi_1(\alpha)}$. Similarly we take 
the chamber $C_{\varphi_0(\alpha)}^-\neq 
C^{\PP^2}_{\varphi_0(\alpha)}$ in $\varphi_0(\alpha)^\perp$ 
sharing the wall 
$\R_{\ge0}\theta(0)_0$ with $C^{\PP^2}_{\varphi_0(\alpha)}$.
We take a real number $0<\e<1$ small enough
such that $\theta(1-\e)_1\in C^{\PP^2}_{\varphi_1(\alpha)}$,
$\theta(1+\e)_1\in C_{\varphi_1(\alpha)}^+$
and $\theta(\e)_0\in C^{\PP^2}_{\varphi_0(\alpha)}$, 
$\theta(-\e)_0\in C_{\varphi_0(\alpha)}^-$. 
\begin{lem}\label{th}
The following hold.\\
$(1)$\ $M_{B}(-\varphi_1(\alpha),
\theta(1+\e)_1)
\neq\emptyset$ for $n\ge 1$.\\
$(2)$\ $M_{B}(-\varphi_0(\alpha),
\theta(-\e)_0)
\neq\emptyset$ for $n\ge 3$.
\end{lem}
\proof
$(1)$ For any $N\in M_{B}(
-\varphi_1(\alpha),\theta(1-\e)_1)$,
we show that the dual vector 
space $N^\ast:=\Hom_\C(N,\C)$  
has a natural $B$-module 
structure and belongs
to $M_{B}(-\varphi_1(\alpha),
\theta(1+\e)_1)$ as follows. We put 
$N^\ast v^{\ast}_i:=\Hom_{\C}(Nv^{\ast}_{2-i},\C)$ 
and define $\gamma^{\ast}_i|_{N^\ast}$ 
and $\delta^{\ast}_j|_{N^\ast}$ by pull 
backs of $\delta^{\ast}_i|_N$ 
and $\gamma^{\ast}_j|_N$, respectively. 
Any surjection $N^{\ast}\to (N')^{\ast}$ 
corresponds to a submodule $N'$ of $N$ and
\begin{equation}\label{conv1}
\dimv((N')^{\ast})=(\dim_\C N'v^{\ast}_2,
\dim_\C N'v^{\ast}_1,
\dim_\C N'v^{\ast}_0). 
\end{equation}
On the other hand,
from (\ref{theta1}) we have 
\begin{equation}\label{conv2}
\theta(1+\e)_1=\e(-2n-1,n,0)+
\frac{n-(n+1)\e}{n}(-n,0,n)
\in\Hom_{\Z}(K(B),\R).
\end{equation}
By (\ref{conv1}) and (\ref{conv2}),
we have the following equality 
\begin{equation}\label{the}
\theta(1+\e)_1((N')^{\ast})=
-\left(\e\theta(0)_1
+\frac{n-(n+1)\e}{n}\theta(1)_1
\right)(N').
\end{equation} 
Since by Lemma~\ref{C_0},
we see that $\theta(1-\e)_1$ and $\e\theta(0)_1
+\frac{n-(n+1)\e}{n}\theta(1)_1$ belong to 
the same chamber $C^{\PP^2}_{\varphi_1(\alpha)}$ for 
$\e$ small enough,
the right hand side of (\ref{the})
is non-positive for any submodule $N'$
of $N\in M_{B}(-\varphi_1(\alpha),\theta(1-\e)_1)$.
We have $\theta(1+\e)_1((N')^{\ast})\le 0$ for any
surjection $N^\ast\to(N')^\ast$. 
Thus $N^{\ast}$ 
belongs to $M_{B}(-\varphi_1(\alpha),\theta(1+\e)_1)$.\\
$(2)$ For $n\ge 3$ we take an element 
$\mathcal{I}_Z\in(\PP^2)^{[n]}$
such that Supp$(\mo_{\PP^2}/\mathcal{I}_Z)$
is not contained in any line $\ell$
on $\PP^2$. Hence we have $\Hom_{\PP^2}
(\mo_{\PP^2},\mathcal{I}_Z(1))=0$.
Below we show that this implies that the $B$-module 
$M:=\Phi_0(\mathcal{I}_Z(1)[1])\in M_{B}
(-\varphi_0(\alpha),\theta(\e)_0)$ is also
$\theta(-\e)_0$-semistable.
For any $B$-submodule $M'\subset
M$, if $\theta(0)_0(M')>0$ then 
by taking $\e$ small enough we have
$\theta(-\e)_0(M')>0$ and $M'$ does not break
$\theta(-\e)_0$-semistability of $M$.
If $\theta(0)_0(M')
=0$, then from (\ref{theta2})
$\dimv M'=(n,\ast,n-1)$ or $(0,\ast,0)$.
However the latter case contradicts
the fact that $\Hom_{B}(\C v_1,M)\cong
\Hom_{\PP^2}(\mo_{\PP^2},\mathcal{I}_Z
(1))=0$. Hence we have $\dimv M'=
(n,l,n-1)$ with $0\le l\le 2n$ and 
$\theta(-\e)_0
(M')\ge 0$. Thus $M$ is 
$\theta(-\e)_0$-semistable. 
\endproof

For $\theta_k\in C_{\varphi_k(\alpha)}^{\PP^2}$,
we have natural morphisms
\begin{equation}\label{mor}
(\PP^2)^{[n]}\cong M_{B}
(-\varphi_k(\alpha),\theta_k) 
\to M_{B}(-\varphi_k(\alpha),
\theta(k)_k)
\end{equation}
for $k=0,1$, since $\R_{\ge0}\theta(1)_1$ 
and $\R_{\ge0}\theta(0)_0$ are 
walls of 
the chamber $C^{\PP^2}_{\varphi_1(\alpha)}$ 
and $C^{\PP^2}_{\varphi_0(\alpha)}$, respectively.
We study the Stein 
factorization $\pi'_k\colon 
(\PP^2)^{[n]}\to Y_k$ of the 
above morphism (\ref{mor}) 
for each $k=0,1$.
Since by Lemma~\ref{th}, for 
$n\ge 3$ our situations 
satisfy the assumptions in 
\cite[Theorem~(3.3)]{Th}, 
we see that $\pi'_1$ and 
$\pi'_0$ are birational 
morphisms and have the 
following diagram: 
\begin{equation}
\label{diagram}
\xymatrix{
M_{B}(-\varphi_0(\alpha),
\theta(-\e)_0)\ar[dr]  & &
\ar@{-->}[ll]_{\kappa }
(\PP^2)^{[n]}\ar[dl] 
\ar[dl]_{\pi'_0} 
\ar[dr]^{\pi'_1}& \\
&  Y_0& & Y_1.}
\end{equation}

\begin{thm}
\label{wall-crossing}
The following 
hold.\\
$(1)$ There exists 
an isomorphism 
$Y_1\cong S^{n}(\PP^2)$
and via this isomorphism, 
the morphism $\pi'_1$ coincide 
with the Hilbert-Chow
morphism $\pi$.\\
$(2)$ For $n\ge 3$, the morphism $\pi'_0$ is 
the contraction morphism 
of the extremal ray $\R_{>0}
\zeta_\ell$. Hence $\pi'_0$ 
coincide with $\psi$ defined in 
\S~$\ref{eff}$ up to isomorphism.
\end{thm}
\proof
$(1)$ We take two elements $\mathcal{I}_Z,
\mathcal{I}_{Z'}\in 
(\PP^2)^{[n]}$. 
We show that if 
$\text{Supp}(Z)=\text{Supp}(Z')$, then 
$\Phi_1(\mathcal{I}_Z(1)[1])$ and 
$\Phi_1(\mathcal{I}_{Z'}(1)[1])$ are 
S-equivalent  
$\theta(1)_1$-semistable 
$B$-modules. By Proposition~\ref{S-eq}
this implies that $\pi'_1$ contracts 
the curve $\beta_n$ to one point. This
shows that the morphism 
$\pi'_1$ coincides with the Hilbert-Chow 
morphism $\pi$ via an isomorphism 
$Y_1\cong S^{n}(\PP^2)$, since 
the Picard number of $(\PP^2)^{[n]}$ 
is two ($n\ge 2$).

We put 
$\text{Supp}(\mo_Z)=
\text{Supp}(\mo_{Z'})=\{x_1,\ldots,x_n\}$
and consider a filtration 
of $\mathcal{I}_Z(1)[1]$ in $\A_1$. 
We put $Z_0:=Z\in(\PP^2)^{[n]}$ and 
inductively define $Z_{i+1}\in(\PP^2
)^{[n-i-1]}$ from $Z_i$ by the following
exact sequence in $\Coh(\PP^2)$
\begin{equation}\label{fil}
0\to \mo_{Z_{i+1}}\to \mo_{Z_i}\to
\mo_{x_{i+1}}\to0
\end{equation}
for $i=0,\ldots, n-2$.
We have $\mo_{Z_{n-1}}=\mo_{x_n}$ and 
$\mo_{x_i}\in\A_1$ for any $i$
by (\ref{skys}). 
By (\ref{fil}) we have 
$\mo_{Z_i}\in\A_1$ for $i=0,
\ldots,n-1$. Hence (\ref{fil}) is also
exact in $\A_1$. On the other hand,
from the exact sequence in $\Coh(\PP^2)$
\begin{equation}\label{prefil}
0\to\mathcal{I}_{Z}\to\mo_{\PP^2}
\to\mo_Z\to0
\end{equation}
we have an exact sequence in $\A_1$
\begin{equation}\label{fil1}
0\to\mo_{Z}\to\mathcal{I}_Z(1)[1]
\to\mo_{\PP^2}(1)[1]\to 0.
\end{equation}

Since $\dimv(\Phi_1(\mo_{\PP^2}(1)[1]))=
(0,1,0)$ and $\dimv(\Phi_1(\mo_x))=
(1,2,1)$ for any closed point $x
\in\PP^2$, we have $\theta(1)_1(
\Phi_1(\mo_{\PP^2}(1)[1]))=\theta(1)_1
(\Phi_1(\mo_x))=0$ from (\ref{theta1}). 
Furthermore from (\ref{fil})
we have $\theta(1)_1(\Phi_1(\mo_{Z_i}))=0$
for any $i$. 
Hence (\ref{fil}) and (\ref{fil1}) give 
a Jordan-H\"{o}lder filtration of 
$\Phi_1(\mathcal{I}_Z(1)[1])$ with 
$\theta(1)_1$-stable quotients 
$\{\Phi_1(\mo_{\PP^2}(1)[1]), \Phi_1(\mo_{x_1}),\ldots,
\Phi_1(\mo_{x_n}) \}$. 
This set only depends on $\text{Supp}
(Z)$. 
Thus
$\Phi_1(\mathcal{I}_{Z}(1)[1])$ and 
$\Phi_1(\mathcal{I}_{Z'}(1)[1])$
represent the same S-equivalence
class of $\theta(1)_1$-semistable $B$-modules.
\\
$(2)$ For a line $\ell$, we take 
an element $\mathcal{I}_Z$
of $\zeta_\ell$. 
As in Lemma~\ref{C_0}, we get an 
exact sequence in $\A_0$
$$0\to\mo_{\PP^2}[1]\to\mathcal{I}_Z(1)[1]\to
\mo_{\ell}(-n+1)[1]\to0$$ 
and $\theta(0)_0(\Phi_0(\mo_{\PP^2}[1]))=\theta(0)_0
(\Phi_0(\mo_{\ell}(-n+1)[1]))=0$. 
Hence by a similar argument as in the 
proof of (1), we see that $\pi'_0$ contracts the curve 
$\zeta_\ell$ on $(\PP^2)^{[n]}$
to one point. 
\endproof

If $n\ge4$, the morphism $\psi$ is 
small and induces a flip in 
the sense of \cite{Th}. 
For general $r>0$ it will be 
shown in \cite{O} that 
$\kappa$ in the above 
diagram~(\ref{diagram})
is the Mori 
flip for $n\gg 0$ and described
by stratified Grassmann bundles.

\noindent\emph{
\bf{Acknowledgement}}\\
The author is grateful to 
his adviser Takao Fujita for 
many valuable comments and 
encouragement.
He thanks Tom Bridgeland, 
Akira Ishii, Emanuele Macri, 
Hiraku Nakajima, Kentaro Nagao, 
Yukinobu Toda, Hokuto Uehara, 
K\={o}ta Yoshioka for valuable 
comments. I wish to thank the 
referee for very careful readings
of the paper and suggesting many
corrections to make the paper more
readable.
This research was supported 
in part by JSPS Global COE program
``Computationism As a Foundation 
of the Sciences''.

Department of Mathematics, Tokyo Institute of Technology, 2-12-1 Oh-okayama, Meguro-ku, 
Tokyo 152-8551, Japan\\
ookawa@math.titech.ac.jp

\begin{thebibliography}{LQZ}
\bibitem[ABL]{ABL}D. Arcara, A. Bertram, M. Lieblich, \emph{Bridgeland-Stable Moduli Spaces for K-Trivial Surfaces}, math.AG/0708.2247.

\bibitem[BBD]{BBD}A. Beilinson, J. Bernstein, P. Deligne, \emph{Faisceaux Pervers},
Ast\'{e}risque 100, Soc. Math de France (1983).

\bibitem[Bo]{Bo}A. I. Bondal, \emph{Representation of associative algebras and coherent sheaves}, Izv. Akad. Nauk
SSSR Ser. Mat. 53 (1989), 25-44; English transl. in Math. USSR-Izv. 34 (1990), no. 1, 23-44.

\bibitem[Br1]{B2}T. Bridgeland, \emph{Stability conditions on triangulated categories}, Ann. of Math. (2), 166 (2007), no. 2, 317-345, also
arXiv:math/0212237.

\bibitem[Br2]{B3}T. Bridgeland, \emph{Stability conditions on K3 surfaces}, Duke Math. J., 141 (2008), no. 2, 241-291, also
arXiv:math/0307164.

\bibitem[Br3]{B4}T. Bridgeland, \emph{T-structures on some local Calabi-Yau varieties}, 
J. Algebra., 289 (2005), no. 2, 453-483, also
arXiv:math/0502050.



\bibitem[Hu]{H}D. Huybrechts, \emph{Fourier-Mukai Transforms in Algebraic Geometry}, Oxford Mathematical Monographs, Oxford Univ. Press, 2006.

\bibitem[HL]{HL}D. Huybrechts, M. Lehn, \emph{Geometry of moduli spaces of sheaves}, Vol. E31 of Aspects in
Mathematics., Vieweg, 1997.

\bibitem[J]{J}D. Joyce, \emph{Configurations in abelian categories. III. 
Stability conditions and identities}, Adv. Math. 215 (2007), 153-219, also arXiv:math/0410267.

\bibitem[K]{Ki1}A. D. King, \emph{Moduli of representations of finite-dimensional algebras}, Quart.
J. Math. Oxford, 45 (1994), 515-530.

\bibitem[KW]{KW}A. D. King, C. H. Walter, \emph{On Chow rings of fine moduli spaces
of modules}, J. Reine Angew. Math. 461 (1995), 179-187.

\bibitem[LQZ]{LQZ}W. Li, Z. Qin, Q. Zhang, \emph{Curves in the Hilbert schemes of points on surfaces}.  Vector bundles and representation theory (Columbia, MO, 2002),  89-96, Contemp. Math., 322, Amer. Math. Soc., Providence, RI, 2003,
also arXiv:math/0105213.

\bibitem[MW]{MW}K. Matsuki, R. Wentworth, \emph{Mumford-Thaddeus principle on the moduli space of vector bundles on an algebraic surface}, Internat. J. Math., 8 (1997), 97-148, also arXiv:alg-geom/9410016.




\bibitem[O]{O}R. Ohkawa, \emph{Flips of moduli of stable
torsion free sheaves with $c_1=1$ on $\PP^2$}, in preparation.

\bibitem[P]{P}J. Le Potier, \emph{A propos de la construction de l'espace de modules des faisceaux
semi-stables sur le plan projectif}, Bull. Soc. Math. France, 122 (1994), 363-369. 

\bibitem[P2]{P2}J. Le Potier, \emph{Lectures on Vector Bundles}, Translated by A. Maciocia, Cambridge Stud. Adv. Math., 54, Cambridge Univ. Press, Cambridge, 1997.

\bibitem[S]{S}C. Simpson, \emph{Moduli of representations of the fundamental group of a smooth
projective variety. I}, Inst. Hautes \'{E}tudes Sci. Publ. Math. 79 (1994), 47-129. 

\bibitem[Th]{Th}M. Thaddeus, \emph{Geometric invariant theory and flips}, J. Amer. Math. Soc., 9 (1996), no. 3, 691-723,
also arXiv:alg-geom/9405004.

\bibitem[To]{To}Y. Toda, \emph{Moduli stacks and invariants of semistable objects on K3 surfaces}, 
Adv. Math. 217 (2008), no. 6, 2736-2781, also arXiv:math/0703590.

\bibitem[TU]{TU}Y. Toda, H. Uehara, \emph{Tilting generators via ample line bundles}, 
Adv. Math. 223 (2010), no. 1, 1-29, also arXiv:0804.4256.

\end{thebibliography}
\end{document}